\documentclass{article}
\usepackage{amsmath,amssymb}
\newtheorem{Theorem}{Theorem}[section]
\newtheorem{Proposition}{Proposition}[section]
\newtheorem{Lemma}{Lemma}[section]
\newtheorem{Corollary}{Corollary}[section]

\newcommand{\Div}{{\rm div}_x}
\newcommand{\dx}{{\rm d}x}
\newcommand{\dt}{{\rm d}t}
\newcommand{\dxdt}{\dx \ \dt}

\newcommand{\intO}[1]{ \int_{\Omega} #1 \ {\rm d}x}
\newcommand{\intOd}[1]{ \int_{\Omega_{\delta}} #1 \ {\rm d}x}

\newcommand{\D}{{\cal D}}

\newcommand{\ep}{\varepsilon}

\newcommand{\Ov}[1]{ \overline{ #1 } }

\newcommand{\Grad}{\nabla_x}
\newcommand{\rD}{\varrho_{\delta}}
\newcommand{\sD}[1]{ \left\{ #1 \right\}_{\delta > 0} }
\newcommand{\uD}{\vc{u}_{\delta}}
\newcommand{\tD}{\vartheta_{\delta}}
\newcommand{\bD}{\vc{B}_{\delta}}
\newcommand{\bH}{\vc{H}_{\delta}}
\newcommand{\bTheorem}[1]{\bigskip \begin{Theorem} \label{T#1}}
\newcommand{\eT}{\end{Theorem} \bigskip }
\newcommand{\bProposition}[1]{\bigskip \begin{Proposition} \label{P#1}}
\newcommand{\eP}{\end{Proposition} \bigskip }
\newcommand{\bLemma}[1]{\bigskip \begin{Lemma} \label{L#1}}
\newcommand{\eL}{\end{Lemma} \bigskip }
\newcommand{\bCorollary}[1]{\bigskip \begin{Corollary} \label{C#1}}
\newcommand{\eC}{\end{Corollary} \bigskip }
\newcommand{\bFormula}[1]{\begin{equation} \label{#1}}
\newcommand{\eF}{\end{equation}}
\newcommand{\bRemark}[1]{\bigskip \begin{Remark} \label{R#1}}
\newcommand{\eR}{\end{Remark} \bigskip }
\newcommand{\vc}[1]{ {\bf #1} }
\newcommand{\Curl}{{\bf curl}_x}

\newcommand{\tn}[1]{\mathbb{#1}}

\def\ddk{\dot \Delta_k}

\newcommand{\R}{{\mathbb R}}
\newcommand{\Z}{{\mathbb Z}}
\def\cF{{\mathcal F}}

\date{}

\begin{document}

\title{Weak and strong solutions of equations of compressible magnetohydrodynamics}
\author{Xavier Blanc$^1$, Bernard Ducomet$^2$ \bigskip \\
\footnotesize{$^1$ Univ. Paris Diderot, Sorbonne Paris Cit\'e,} \\
\footnotesize{Laboratoire Jacques-Louis Lions, UMR 7598,} \\ 
\footnotesize{UPMC, CNRS, F-75205 Paris, France} \medskip\\
\footnotesize{$^2$ CEA, DAM, DIF,} \\
\footnotesize{D\'epartement de Physique Th\'eorique et Appliqu\'ee,} \\
\footnotesize{F-91297 Arpajon, France}}

\maketitle

\section{Basic equations of magnetohydrodynamics}
\label{b}

A plasma is a mixture of ions, electrons and neutral particles. At the macroscopic level, beyond the kinetic description necessary 
in rarefied situations, dynamics of dense plasmas is governed by the interaction between the fluid components and the electromagnetic fields.
When dealing with collisional plasmas \cite{KrL} \cite{Ku}, it is often appropriate to use one-fluid approximations
which are much simpler than
a complete kinetic theory and relevant for a number of applications (astrophysics \cite{COGI} \cite{SHOR} \cite{ZZJ} \cite{ZIR}, 
plasma physics \cite{DMP} \cite{DEBE}, electrometallurgy \cite{GLBL} etc...).

Moreover a further approximation is obtained when one studies the
interaction between electromagnetic fields and a electric-conducting fluid with electrical effects negligeable compared to
 magnetic effects. This simplified scheme is called magnetohydrodynamics (MHD).
It unifies classical fluid dynamics and magnetism of continuum media in a coupled way: electromagnetic fields induce currents
 in the moving fluid which
produce forces which in turn modify the electromagnetic fields. If one considers a compressible heat-conducting fluid,
 the equations governing the system will be
the compressible Navier-Stokes-Fourier system coupled to the set of Maxwell's equations through suitable momenta and energy sources.

\subsection{Modelling of the magnetic field}

Recall that the electromagnetic field is governed by the Maxwell system \cite{LL} decomposed into
\begin{itemize}
\item the {\it Amp\`ere's law}
\bFormula{A}
\partial_t \vc{D}+\vc{j} = \Curl \vc{H},
\eF
\item the {\it Coulomb's law}
\bFormula{C}
 \Div \vc{D} = \varrho_c,
\eF
\item the {\it Faraday's law}
\bFormula{F}
\partial_t \vc{B} + \Curl \vc{E} = 0,
\eF
\item the {\it Gauss's law}
\bFormula{G}
\Div \vc{B} = 0,
\eF
\end{itemize}
where $\vc{D}$ is the electric induction, $\vc{B}$ is the magnetic induction, $\vc{E}$ is the electric field, $\vc{H}$ is the magnetic field,
 $\varrho_c$ is the density of electric charges and $\vc{j}$ the electric current.
We suppose that these quantities are linked by the constitutive relations
\bFormula{mueps}
 \vc{B} = \mu  \vc{H}\ \ \ \ \mbox{and}\ \ \ \  \vc{D} = \varepsilon  \vc{E},
\eF
where $\mu$, the magnetic permeability and $\varepsilon$ the dielectric permittivity are positive quantities, generally depending on the electromagnetic field.
Furthermore one observes that the electric density charge $\varrho_c$ and the current density $\vc{j}$ are
 interrelated through the {\it electric charge conservation}
\bFormula{ECC}
\partial_t \varrho_c +\Div  \vc{j}=0,
\eF
however we assume in the following that electric neutrality holds: then $\varrho_c=0$.

Finally we suppose that the magnetic induction vector $\vc{B}$ is related to the
electric field $\vc{E}$ and the macroscopic fluid velocity $\vc{u}$ via {\it Ohm's law}
\bFormula{O}
\vc{j} = \sigma ( \vc{E} + \vc{u} \times \vc{B}),
\eF
where the electrical conductivity $\sigma = \sigma(\varrho, \vartheta, \vc{H} )$ of the fluid is a positive quantity.

Now the {\it magnetohydrodynamic approximation} \cite{C} \cite{KL} \cite{Wo} consists in neglecting the term $\partial_t \vc{D}$ in \eqref{A} and
we obtain
\bFormula{Abis}
 \vc{j} = \Curl \vc{H}.
\eF
This approximation corresponds to a strongly non-relativistic regime and has been justified in the 2-dimensional case by S. Kawashima and Y. Shizuta in \cite{KS},
and more recently in 3-dimensional space 
by S. Jiang and F. Li (see \cite{JL} for isentropic flows and \cite{JLbis} for the general case).

Accordingly, equation \eqref{F} can be written in the form
\bFormula{b7}
\partial_t \vc{B} + \Curl( \vc{B} \times \vc{u} ) + \Curl\left( \frac{1}{\sigma}\ \Curl \left(\frac{1}{\mu}\ \vc{B}\right)\right) = 0,\
\eF
where $\mu = \mu( |\vc{H}| ) > 0$.

Setting
\bFormula{EM}
 {\mathcal M}(s)=\int_0^s \tau\partial_{\tau}\left( \tau\mu(\tau)\right)\ d\tau,
\eF
 equation \eqref{F} becomes
\bFormula{EMEmhd}
 \partial_t {\mathcal M}(|\vc{H}|)
+ \vc{j}\cdot \vc{E}=\Div \left(\vc{H}\times\vc{E}\right).
\eF
Supposing now that we study the previous system in a bounded region $\Omega\subset {\R}^3$, we add boundary conditions
\bFormula{bcBE}
\left. \vc{B}\cdot \vc{n}\right|_{\partial\Omega}  = 0,\ \ \ \ \left. \vc{E}\times \vc{n}\right|_{\partial\Omega}  = \vc{0}.
\eF
and initial conditions
\bFormula{icBE}
\vc{B}(0,x)  = \vc{B}_0(x),\ \ \ \  \vc{E}(0,x)  = \vc{E}_0(x)\ \ \ \mbox{for any}\ x\in\Omega.
\eF
In accordance with the basic principles of continuum mechanics, the magnetofluid dynamics is described (see Cabannes \cite{C})
by the {\it Navier-Stokes-Fourier system} of equations with magnetic sources. This includes
{\it mass conservation}
\bFormula{b8}
\partial_t \varrho + \Div (\varrho \vc{u}) = 0,
\eF
 the {\it momentum balance}
 \begin{equation}
   \label{b9.0}
\partial_t (\varrho \vc{u}) + \Div (\varrho \vc{u} \otimes \vc{u}) = \Div \tn{T} +\varrho_c \vc{E}+ \vc{j} \times
\vc{B} + \vc{f},   
 \end{equation}
and the {\it energy  balance}
\bFormula{b18bis}
\partial_t \Big( \frac{1}{2} \varrho |\vc{u}|^2 + \varrho e  \Big) 
+\Div \Big( (\frac{1}{2} \varrho |\vc{u}|^2 + \varrho e + p) \vc{u}  - \tn{S} \vc{u} \Big)+ \Div \vc{q} = \vc{E}\cdot\vc{j} + g,
\eF
where $e$ is the specific internal energy, $\vc{q}$ denotes the heat flux, $\tn{T}$ is the Cauchy stress tensor,
 $\vc{f}$ is a given body force, $g$ is a given energy source 
and the Lorentz force $\varrho_c\vc{E}+\vc{j} \times \vc{B}$ is imposed by the electromagnetic field. In the
magnetohydrodynamic approximation, it may be shown \cite{C} that this approximation also implies that
$\varrho_c \vc{E}$ is negligible compared to ${\mathbf j}\times {\mathbf B}$. Therefore, the momentum balance \eqref{b9.0} reads
 \begin{equation}
   \label{b9}
\partial_t (\varrho \vc{u}) + \Div (\varrho \vc{u} \otimes \vc{u}) = \Div \tn{T} + \vc{j} \times
\vc{B} + \vc{f}.   
 \end{equation}
Furthermore the Cauchy stress tensor $\tn{T}$
is given by {\it Stokes' law}
\bFormula{b12}
\tn{T} = \tn{S} - p \tn{I},
\eF
where $p$ is the pressure, and the symbol $\tn{S}$ stands for the viscous stress tensor, given by
\bFormula{b13}
\tn{S} = \nu \left( \Grad \vc{u} + \Grad \vc{u}^t - \frac{2}{3} \Div (\vc{u}) \ \tn{I} \right) + \eta \ \Div (\vc{u}) \ \tn{I},
\eF
with the shear viscosity coefficient $\nu$, and the bulk viscosity coefficient $\eta$.

Finally we impose a no-slip boundary condition for the velocity
\bFormula{ns}
\left.\vc{u}\right|_{\partial\Omega} = 0.
\eF
Using \eqref{EMEmhd} and \eqref{b18bis} the energy of the system satisfies
\bFormula{b14}
\partial_t \Big( \frac{1}{2} \varrho |\vc{u}|^2 + \varrho e + {\mathcal M}(|\vc{H}|)\Big) +
\Div \left( \Big(\frac{1}{2} \varrho |\vc{u}|^2 + \varrho e + p\Big) \vc{u} + \vc{E} \times \vc{B} - \tn{S} \vc{u}+\vc{q} \right) = 0,
\eF
where $e$ is the specific internal energy, and $\vc{q}$ denotes the heat flux.

We prescribe a no-flux boundary condition for $\vc{q}$
\bFormula{nf}
\vc{q} \cdot \vc{n}|_{\partial \Omega} = 0.
\eF
Using the previous boundary conditions we deduce that the (conserved) total energy of the system satisfies
\bFormula{b15}
{{\rm d} \over {\rm d}t} \int_{\Omega}\Big( {1 \over 2} \varrho |\vc{u}|^2 + \varrho e + {\mathcal M}(|\vc{H}|) \Big)\ \dx = 0.
\eF
In the variational formulation it will be more convenient to replace
\eqref{b14} by the {\it entropy balance}
\bFormula{b16}
\partial_t (\varrho s) + \Div (\varrho s \vc{u} ) + \Div \Big(
{\vc{q} \over \vartheta } \Big) = r,
\eF
with the specific entropy $s$, and the entropy production rate $r$, for which the second law of thermodynamics requires $r \geq 0$.

The specific entropy $s$, the specific internal energy $e$, and the pressure are interrelated through the
second principle of thermodynamics
\bFormula{b17}
\vartheta D s = D e + D \Big( {1 \over \varrho} \Big) p,
\eF
where $D$ stand for the differential with respect to the variable $\varrho$, $\vartheta$, which implies some compatibility conditions 
between $e$ and $p$ ( Maxwell's relation).

If the motion is smooth, the sum of the kinetic and magnetic energy satisfies
\bFormula{b18}
\partial_t \left(\frac{1}{2}  \varrho |\vc{u}|^2 + {\mathcal M}(|\vc{H}|)\right)
+\Div \left(  \Big(\frac{1}{2} \varrho |\vc{u}|^2 +  p\Big) \vc{u} + \vc{E} \times \vc{B} - \tn{S} \vc{u} \right)
 =p \ \Div (\vc{u}) - \tn{S} : \Grad \vc{u} - \frac 1 {\mu\sigma} | \Curl \vc{B} |^2.
\eF
 Indeed using \eqref{b17} one checks that
\bFormula{b19}
r = \frac{1}{\vartheta}
\left( \tn{S} : \Grad \vc{u} + \frac{1}{\sigma} |\Curl \vc{H}|^2 - \frac{\vc{q} \cdot \Grad \vartheta}{ \vartheta} \right).
\eF
If the motion is not (known to be) smooth, validity of \eqref{b18bis} is no longer guaranteed and the dissipation rate of mechanical 
energy may exceed the value $\tn{S} : \Grad \vc{u}$.

Accordingly, we replace equation \eqref{b16} by inequality
\bFormula{b20}
\partial_t (\varrho s) + \Div (\varrho s \vc{u} ) + \Div \left(
\frac{\vc{q}}{ \vartheta } \right) \geq
\frac{1}{\vartheta}
\Big( \tn{S} : \Grad \vc{u} + \frac{1}{\sigma} |\Curl \vc{H}|^2 - \frac{\vc{q} \cdot \Grad \vartheta}{\vartheta} \Big),
\eF
to be satisfied together with the total energy balance \eqref{b15}.

\subsection{Hypotheses on constitutive relations}

In order to close the system one needs first to postulate a constitutive equation relating the pressure to the other variables. 
Several possibilities have been considered. We choose here, in the spirit of \cite{DUFE2} (however see \cite{FY} \cite{EF70} \cite{HW1} for other possible choices) to focus on high-temperature phenomena appearing
 in various astrophysical contexts (stars, planetar magnetospheres, neighborhood of black holes etc...) 
 \cite{COGI} \cite{SHOR} \cite{ZZJ} \cite{ZIR}, or plasma physics applications (tokamaks, inertial confinement fusion devices etc...) \cite{KrL} \cite{Ku}.

In that case, the so called equilibrium diffusion limit of radiation hydrodynamics involves a modified pressure law where
the pressure $p_F(\varrho, \vartheta)$ in the fluid is augmented by a radiation
component $p_R( \vartheta)$ related to the absolute temperature $\vartheta$ through the Stefan-Boltzmann law
\bFormula{b1}
p_R( \vartheta) = \frac{a}{3} \vartheta^4,\ \mbox{with a constant}\ a > 0
\eF
(see, for instance, Chapter 15 in \cite{EGH}). Accordingly, the specific internal energy of the fluid must be supplemented by
a term
\bFormula{b2}
e_R( \vartheta) = e_R(\varrho, \vartheta) = \frac{a}{\varrho} \vartheta^4,
\eF
where $\varrho$ is the fluid density and similarly, the heat conductivity of the fluid is enhanced by a radiation component 
(see \cite{BUDE} \cite{MIMI})
\bFormula{b3}
\vc{q}_R( \vartheta) = - \kappa_R \vartheta^3 \Grad \vartheta,\ \mbox{with a constant}\  \kappa_R > 0
\eF
Finally the total pressure $p$ takes the form
\bFormula{b21}
p(\varrho, \vartheta) = p_F(\varrho, \vartheta) + p_R( \vartheta),
\eF
where the radiation component is given by \eqref{b1}. Similarly, the specific internal energy
reads
\bFormula{b22}
e (\varrho, \vartheta) = e_F (\varrho, \vartheta) + e_R(\varrho, \vartheta),
\eF
with $e_R$ determined by \eqref{b2}.

Furthermore we suppose that the gas is perfect and monoatomic
\bFormula{b23}
p_F(\varrho, \vartheta) = \frac{2}{3} \varrho e_F(\varrho, \vartheta).
\eF
The following hypothesis expresses the convexity of free energy
\bFormula{b24}
\frac{\partial p_F (\varrho, \vartheta)}{\partial \varrho} > 0,\
\frac{\partial e_F (\varrho, \vartheta)}{\partial \vartheta} > 0 \ \mbox{for all}\ \varrho, \vartheta > 0,
\eF
and finally, we suppose
\bFormula{b25}
\liminf_{\varrho \to 0+}  p_F (\varrho , \vartheta)=0,\ \
\liminf_{\varrho \to 0+} \frac{\partial p_F (\varrho , \vartheta)}{\partial \varrho } > 0 \ \mbox{for any fixed}\ \vartheta > 0,
\eF
and
\bFormula{b26}
\liminf_{\vartheta \to 0+}  e_F (\varrho , \vartheta)=0,\ \ 
\liminf_{\vartheta \to 0+} \frac{\partial e_F (\varrho, \vartheta)}{\partial \varrho} > 0 \ \mbox{for any fixed}\ \varrho > 0.
\eF
As one can check by solving \eqref{b17} with \eqref{b23}, there exists a function 
$P_F \in C^1(0, \infty)$ such that 
\bFormula{p1}
p_F (\varrho, \vartheta) = \vartheta^{\frac{5}{2}} P_F \Big( \frac{\varrho}{\vartheta^{\frac{3}{2}}} \Big).
\eF
Consequently
\[
\frac{\partial e_F (\varrho, \vartheta)}{\partial \vartheta} = \frac{3}{2} \frac{1}{Y} \Big( \frac{5}{3} P_F(Y) - P'_F(Y) Y \Big) ,\ 
 \mbox{with}\ \ Y = \frac{\varrho}{\vartheta^{\frac{3}{2}}},
\]
so according to \eqref{b24} 
\bFormula{p2}
P_F'(z) > 0,\ \frac{5}{3} P_F(z) - P_F'(z) z > 0 \ \mbox{for any}\ z > 0,
\eF
where the latter inequality yields
\bFormula{p3}
\Big( \frac{ P_F(z)}{z^{\frac{5}{3}}} \Big)' < 0 \ \mbox{for all}\ z > 0. 
\eF
In particular
\bFormula{p4}
\frac{ P_F (z)}{z^{\frac{5}{3}}} \to p_{\infty} > 0 \ \mbox{for}\ z \to \infty.
\eF
According to \eqref{b25}, \eqref{b26}, we see that $P_F \in C^1[0, \infty)$. Moreover
\bFormula{p5}
0 < \liminf_{z \to 0+} \frac{1}{z} \Big( \frac{5}{3}P_F(z) - P'_F(z) z \Big) \leq 
\limsup_{z \to 0+} \frac{1}{ z} \Big(\frac{5}{3} P_F(z) - P'_F(z) z \Big) < \infty,
\eF
\bFormula{p6}
\limsup_{z \to \infty}\frac{1}{z} \Big(\frac{5}{3} P_F(z) - P'_F(z) z \Big) < \infty,
\eF
and
\bFormula{p7}
\lim_{z \to 0+} P_F(z) = 0,\  \lim_{z \to 0+} P'_F(z) > 0,\ 
\lim_{z \to \infty} {P'_F(z) \over z^{\frac{2}{3}} } = \frac{5}{3}p_{\infty} > 0.
\eF
Now it follows from \eqref{p7} that 
\[
P'_F(z) \geq c z^{\frac{2}{3}} \ \mbox{for all}\ z > 0, 
\ \mbox{and a certain}\ c > 0.
\]
Therefore there exists $p_c > 0$ such that the mapping
\bFormula{p8}
\varrho \mapsto 
p_F(\varrho, \vartheta) - p_c \varrho^{\frac{5}{3}}, 
\eF
is a non-decreasing function of $\varrho$ for any fixed $\vartheta > 0$.

Accordingly, using equation \eqref{b17} we have
\bFormula{p9}
s(\varrho, \vartheta) = s_F(\varrho, \vartheta) + s_R(\varrho, \vartheta),
\eF
where
\bFormula{p10}
s_F (\varrho, \vartheta) = S_F \Big( \frac{\varrho}{\vartheta^{\frac{3}{2}}} \Big),\ \mbox{with}
\ S_F'(z) = -\frac{3}{2} {\frac{5}{3} P_F(z) - P'_F(z) z \over z^2 },
\eF
and
\bFormula{p11}
\varrho s_R (\varrho, \vartheta) = {4 \over 3} a \vartheta^3.
\eF
Note that after \eqref{p2} and \eqref{p5}, $S_F$ is a decreasing function such that
\[
\lim_{z \to 0+} z S'_F(z) = - {2 \over 5} P'_F(0) < 0, 
\]
whence, normalizing by $S_F(1) = 0$, we get 
\bFormula{p12}
-c_1 \log(z) \leq S_F(z) \leq - c_2 \log(z),\ c_1 > 0, 
\ \mbox{for all}\ 0 < z \leq 1,
\eF
and
\bFormula{p13}
0 \geq S_F(z) \geq - c_3 \log(z) \ \mbox{for all}\ z \geq 1.
\eF
Using the second law of thermodynamics, we see that viscosity coefficients $\nu$ and $\eta$ are non-negative quantities.

In addition, we shall suppose that
\bFormula{b27}
\nu = \nu( \vartheta, |\vc{H}| ) > 0,\ \eta = \eta (\vartheta, | \vc{H} |) \geq 0,
\eF
where the bulk viscosity $\nu$ satisfies some technical, but physically relevant, coercivity conditions to be specified below.

Similarly, the heat flux $\vc{q}$ obeys Fourier's law
\bFormula{b28}
\vc{q} = \vc{q}_R + \vc{q}_F,
\eF
where the radiation heat flux $\vc{q}_R$ is given by \eqref{b3} and
\bFormula{b29}
\vc{q}_F = - \kappa_F(\varrho, \vartheta, |\vc{H}|) \Grad \vartheta,\ \kappa > 0.
\eF
For the sake of simplicity we assume the transport coefficients $\nu$, $\kappa_F$, $\mu$ and $\sigma$ to admit a common temperature scaling,
namely
\bFormula{b31}
c_1 (1 + \vartheta^{\alpha}) \leq
\kappa_F(\varrho, \vartheta, |\vc{H}|),\ \nu(\vartheta, |\vc{H}|), \sigma^{-1} (\varrho, \vartheta, |\vc{H}| ) \leq
c_2 (1 + \vartheta^{\alpha}),\ c_1,\ c_2 > 0,
\eF
with $\alpha \geq  1$ to be specified below.
Similarly, we suppose that
\bFormula{b32}
0 \leq \eta(\vartheta, |\vc{H}|) \leq c_3 (1 + \vartheta^{\alpha}),
\eF 
and that the magnetic permeability has the symbol property
\bFormula{b32bis}
\underline{c}_k s(1 + s)^{-k} \leq\partial_s^k(s\mu( s)) \leq \overline{c}_k s(1 + s)^{-k},
\eF
for any $s\geq 0$ and for $k=0,1$, with $\underline{c}_k,\ \overline{c}_k>0$.

We finally suppose that external sources are absent: $\vc{f}=0$ and $g=0$.

\vskip0.5cm
Let us conclude this introduction by briefly fixing the limits of our study.
As our framework is explicitely viscous and compressible, all the important works about the incompressible model are beyond our scope.
However the reader may consult, among standard references: Ladyzhenskaya and Solonnikov \cite{LASO}, Duvaut and Lions \cite{DULI1} \cite{DULI},
 Sermange and Temam \cite{SETE} and Cannone et al. \cite{CA}
for reviews of the incompressible setting and Germain and Masmoudi \cite{GEMA} (with references therein) for the inviscid case.
Finally we refer the reader to Eringen and Maugin \cite{EM} for more complex situations coupling electrodynamics to continua.
\vskip0.5cm
The article is structured as follows: in Section \ref{weak} we consider the main existence result of weak solutions,
in Section \ref{other} we present related models (barotropic or polytropic) and in Section \ref{singular} we introduce asymptotic models. In Section \ref{strong}
we present the main results concerning existence and uniqueness of strong solutions including critical aspects and finally results for one dimensional models are
given in the last Section \ref{1D}.

\section{Weak solutions}
\label{weak}

\subsection{Variational formulation}
\label{v}

As in the pure hydrodynamical model \cite{DUFE1} the approach is based on the concept of variational solutions.

1- The mass conservation is expressed through the weak form
\bFormula{v1}
\int_0^T \intO{ \Big( \varrho \partial_t \varphi + \varrho \vc{u} \cdot \Grad \varphi  \Big) } \ \dt +
\intO{ \varrho_0 \varphi(0, \cdot) }  = 0,
\eF
to be satisfied for any  test function $\varphi \in \D([0,T) \times \R^3)$ where $\varrho_0$ is the initial density.

Assuming that both $\varrho$ and the flux $\varrho \vc{u}$ are locally integrable on $[0,T) \times \Ov{\Omega}$ equation \eqref{v1}
makes sense. Moreover if $\varrho$ belongs to $L^{\infty}(0,T; L^{\gamma}(\Omega))$ for a $\gamma > 1$, then \cite{NOST4}
$\varrho \in C([0,T]; L^{\gamma}_{weak}(\Omega))$ i.e.
\[
\varrho(t) \to \varrho_0 \ \mbox{weakly in}\ L^{\gamma}(\Omega) \ \mbox{as}\ t \to 0,
\]
and the total mass 
\bFormula{v2}
M_0 = \intO{ \varrho (t) } = \intO{ \varrho_0 },
\eF
is a constant of motion.
  
For $\gamma \geq 2$ one can use the DiPerna-Lions theory \cite{DL} to show 
that the pair $(\varrho, \vc{u})$ satisfies the renormalized equation 
\bFormula{v3}
\int_0^T \intO{ \Big( b(\varrho) \partial_t \varphi + b(\varrho) \vc{u} \cdot \Grad \varphi + 
(b(\varrho) - b'(\varrho) \varrho) \Div (\vc{u}) \varphi \Big) } \ \dt + \intO{ 
b(\varrho_0) \varphi (0, \cdot) } = 0,
\eF
for any $\varphi \in \D([0,T) \times \R^3)$, and for any continuously differentiable function $b$ whose derivative vanishes for large arguments.
It can also be shown that any renormalized solution $\varrho$ belongs to the space $C([0,T]; L^1(\Omega))$ (see \cite{DL}).

2- The variational formulation of the momentum equation reads now
\[
\int_0^T \intO{ \Big( (\varrho \vc{u}) \cdot \partial_t \varphi +(\varrho \vc{u} \otimes \vc{u}): \Grad \varphi + p \ \Div (\varphi) \Big)}\ \dt
\]
\bFormula{v4}
 =\int_0^T \intO{ \Big( \tn{S} : \Grad \varphi - (\vc{j} \times \vc{B}) \cdot \varphi \Big) } \ \dt - \intO{ (\varrho \vc{u})_0\cdot\varphi(0,\cdot)},
\eF
which must be satisfied for any vector field $\varphi \in \D([0,T) \times \Omega; \R^3)$. Once more we suppose that the quantities 
$\varrho \vc{u}$, $\varrho \vc{u} \otimes \vc{u}$, $p$, $\tn{S}$, $\varrho \Grad \Psi$, and $\vc{J} \times \vc{B}$ are at least locally integrable 
on the set $[0,T) \times \Omega$.

The initial value of the momentum $(\varrho \vc{u})_0$ must satisfy
\bFormula{v6}
(\varrho \vc{u})_0 = 0 \ \mbox{a.e. on the set}\ \{ \varrho_0 = 0 \} .
\eF
Finally using the no-slip boundary condition \eqref{ns} we get
\bFormula{v7}
\vc{u} \in L^2(0,T; W^{1,2}_0 (\Omega; \R^3)).
\eF
Observe that our choice of the function space has been motivated by the a priori bounds resulting from the entropy balance specified below.

3- According to \eqref{b20} the weak form of the entropy production reads
\bFormula{v8}
\int_0^T \intO{ \Big(  \varrho s \ \partial_t \varphi + \varrho s \vc{u} \cdot \Grad \varphi +
{\vc{q} \over \vartheta} \cdot \Grad \varphi \Big) } \ \dt \leq
\eF
\[
\int_0^T \intO{
 {1 \over \vartheta}  \Big( {\vc{q} \cdot \Grad \vartheta \over \vartheta} -
\tn{S} : \Grad \vc{u} - {1 \over {\mu^2\sigma}} |\Curl B|^2 \Big) \varphi }\ \dt 
- \intO{
(\varrho s)_0 \varphi (0,\cdot )},
\]
for any test function $\varphi \in \D([0,T) \times \R^3)$, $\varphi \geq 0$.

We observe that the presence of $\vartheta$ in the denominator forces this quantity to be positive on a set of full 
measure for \eqref{v8} to make sense. 

4- Finally Maxwell's equation \eqref{b7} is replaced by
\bFormula{v9}
\int_0^T \intO{ \Big( \vc{B} \cdot \partial_t \varphi
 - \left(\vc{B} \times \vc{u} + \frac{1}{\sigma} \Curl \frac{\vc{B}}{\mu}\right)  \cdot \Curl \varphi \Big) } \ \dt +
\intO{ \vc{B}_0  \cdot \varphi (0, \cdot) },
\eF
to be satisfied for any vector field $\varphi \in \D([0,T) \times \R^3)$.

Hereafter using boundary conditions \eqref{bcBE} and \eqref{v7} one must take
\bFormula{v10}
\vc{B}_0 \in L^2(\Omega),\ \ \ \ \Div\left(\vc{B}_0\right) = 0 \ \mbox{in}\ \D'(\Omega),\ \ \ \ \ \vc{B}_0 \cdot \vc{n}|_{\partial \Omega} = 0.
\eF
By virtue of Theorem 1.4 in \cite{TEM}, $\vc{B}_0$ belongs to the closure of solenoidal fields on $\D(\Omega)$ with 
respect to the $L^2-$norm.

Applying \eqref{b15} and \eqref{b20}, we infer that ${\mathcal M}(|\vc{H}|) \in L^{\infty}(0,T; L^1(\Omega; \R^3))$ and
 $\Curl \vc{H} \in L^2 (0,T; L^2(\Omega; \R^3))$. But using \eqref{EM} and \eqref{b32bis} we deduce that $\vc{B} \in L^{\infty}(0,T; L^2(\Omega; \R^3))$, 
and from \eqref{G} and \eqref{bcBE} we see that 
\[
\Div (\vc{B})(t) = 0 \ \mbox{in}\ \D'(\Omega),\ \ \ \vc{B}(t) \cdot \vc{n}|_{\partial \Omega} = 0 \ \mbox{for a.e.}\ t \in (0,T),
\]
so using Theorem 6.1 in \cite{DULI} we conclude that
\bFormula{v11}
\vc{H} \in L^2(0,T; W^{1,2}(\Omega; \R^3)),\ \ \
\Div (\vc{B}) (t) = 0,\ \ \ \vc{B} \cdot \vc{n}|_{\partial \Omega} = 0 \ \mbox{for a.e.}\ t \in (0,T).
\eF

5- According to \eqref{b14} we shall assume the total energy to be a constant of motion
\bFormula{v12}
E(t)= \intO{ \Big( {1 \over 2} \varrho(t) |\vc{u}(t)|^2 + (\varrho e)(t)
+ {\mathcal M}(|\vc{H}(t)|) \Big)  } = E_0 \ \mbox{for a.e.}\ t \in (0,T),
\eF
where
\bFormula{v13}
E_0 = \intO{ \Big( {1 \over 2  \varrho_0} | (\varrho \vc{u})_0|^2 + (\varrho e)_0
+ {\mathcal M}(|\vc{H}_0|) \Big)  }.
\eF
Of course the initial values of the entropy $(\varrho s)_0$ and the internal energy $(\varrho e)_0$ should 
be chosen consistently with the constitutive relations. A natural possibility consists in fixing 
the initial temperature $\vartheta_0$ and setting
\bFormula{v15}
(\varrho s)_0 = \varrho_0 s(\varrho_0, \vartheta_0),\ 
(\varrho e)_0 = \varrho_0 e(\varrho_0, \vartheta_0).
\eF

\subsection{Global existence of weak solutions for large data}
\label{m}

The following result holds
\bTheorem{m1}
Let $\Omega \subset \R^3$ be a bounded domain with boundary of class $C^{2 + \delta}$, $\delta > 0$.

Suppose that the thermodynamic functions 
\[
p=p(\varrho, \vartheta), \ e = e(\varrho, \vartheta),\  
s = s(\varrho, \vartheta)
\]
are interrelated through \eqref{b17}, where $p$, $e$ can be decomposed 
as in \eqref{b21}, \eqref{b22}, with the components $p_F$, $e_F$ satisfying \eqref{b23}. Moreover, let $p_F(\varrho, 
\vartheta)$, $e_F(\varrho, \vartheta)$ be continuously differentiable functions of positive arguments $\varrho$, $\vartheta$ 
satisfying \eqref{b24}-\eqref{b26}. 

Furthermore, we suppose that the transport coefficients 
\[
\nu = \nu(\vartheta, |\vc{H}|),\ 
\eta = \eta(\vartheta, |\vc{H}|),\ \kappa_F = \kappa_F(\varrho, \vartheta, |\vc{H}|), \ \sigma = 
\sigma(\varrho, \vartheta, |\vc{H}|)\ \mbox{and} \ \mu=\mu( |\vc{H}|),
\]
are continuously differentiable functions of their arguments obeying 
\eqref{b27}-\eqref{b32bis}, with
\bFormula{m1}
1 \leq \alpha < {65 \over 27}.
\eF
Finally, let the initial data $\varrho_0$, $(\varrho \vc{u})_0$, $\vartheta_0$, $\vc{B}_0$ be given so that
\bFormula{m2}
\varrho_0 \in L^{5 \over 3}(\Omega),\ 
(\varrho \vc{u})_0 \in L^1(\Omega; \R^3),\ \vartheta_0 \in L^{\infty} (\R^3),\ \vc{B}_0 \in L^2(\Omega; \R^3),
 \eF
\bFormula{m3}
\varrho_0 \geq 0,\ \vartheta_0  > 0, 
\eF
\bFormula{m4}
(\varrho s)_0 = \varrho_0 s(\varrho_0, \vartheta_0),\ 
{1 \over \varrho_0 } |(\varrho \vc{u})_0 |^2 ,\ (\varrho e)_0 = \varrho_0 e(\varrho_0, \vartheta_0) 
\in L^1(\Omega),
\eF
and
\bFormula{m5}
\Div (\vc{B}_0) = 0 \ \mbox{in}\ \D'(\Omega),\ \vc{B}_0 \cdot \vc{n}|_{\partial \Omega} = 0.
\eF
Then problem \eqref{v1}-\eqref{v10} possesses at least one variational solution $\varrho$, $\vc{u}$, 
$\vartheta$, $\vc{B}$ in the sense of Section~\ref{v} on an arbitrary time interval $(0,T)$.
\eT
Observe that hypothesis \eqref{m1} is a technical restriction. However functional dependence on $\vartheta$
was rigorously justified in \cite{DELE} after an asymptotic analysis
of some kinetic-fluid models. The same temperature scaling on all transport coefficients was also suggested in \cite{GIO}.
 
\subsubsection{Difficulties and methods}
\label{param}

The main difficulty of the problem being the lack 
of  a priori estimates, methods of weak convergence based on the theory of the 
parametrized (Young) measures represent the main tool \cite{PED1} to proceed.
Accordingly, we denote by $\Ov{b(z)}$ a weak limit of any sequence of composed functions
$\{ b(z_n) \}_{n \geq 1}$. More precisely, 
\[
\int_{\R^M} \Ov{ b(y, \vc{z}) } \varphi (y) \ {\rm d} y = 
\int_{\R^M} \varphi (y) \Big( \int_{\R^K} b(y, z) \ {\rm d} \Lambda_y (z) \Big) \ {\rm d}y, 
\]
where $\Lambda_y (z)$ is a parametrized measure associated to a sequence $\{ \vc{z}_n \}_{n \geq 1}$ of vector-valued functions, 
$\vc{z}_n: \R^M \to \R^K$ (see Chapter 1 in \cite{PED1}).

Assume there is a sequence of approximate solutions resulting from a suitable regularization process. The starting point 
is to establish the relation
\bFormula{m6}
\Ov{ \left( p(\varrho, \vartheta) - \Big({4 \over 3} \nu + \eta\Big) \Div (\vc{u}) \right) b(\varrho)} = \left(
\Ov{ p(\varrho, \vartheta) } \ \Ov{ b(\varrho)} - 
\Big({4 \over 3} \Ov{ \nu(\vartheta, |\vc{B}|) } + \Ov{\eta (\vartheta, |\vc{B}|)}\Big) \Div (\vc{u}) \right) \Ov{b(\varrho)},
\eF
for any bounded function $b$. 

The quantity $p - (4/3 \nu + \eta)\Div (\vc{u})$ is called
the effective viscous pressure and relation \eqref{m6} was proved by P.-L. Lions \cite{LI4} 
for the barotropic Navier-Stokes system, with $p = p(\varrho)$, when viscosity coefficients $\nu$
and $\eta$ are constant. The same result for general temperature dependent viscosity coefficients was obtained in 
\cite{EF71} with the help of certain commutator estimates in the spirit of \cite{COME}.

In the present setting, this approach has to be modified in order to accommodate 
the dependence of $\nu$ and $\eta$ on the magnetic field $\vc{B}$.

The propagation of density oscillations can be described by the 
renormalized continuity equation \cite{DL}
\bFormula{m7}
\partial_t b(\varrho) + \Div (b(\varrho) \vc{u}) + \Big( b'(\varrho) \varrho - b(\varrho) \Big) \ \Div (\vc{u}) = 0,
\eF
and its ``weak'' counterpart
\bFormula{m8}
\partial_t \Ov{b(\varrho)} + \Div ( \Ov{b(\varrho)} \vc{u}) + \Ov{ \Big( b'(\varrho) \varrho - b(\varrho) \Big) \ \Div (\vc{u})} = 0.
\eF
To prove \eqref{m7} one must show first boundedness of the oscillations defect measure
\bFormula{m9}
{\bf osc}_{\gamma + 1}[\varrho_n \to \varrho](Q) = \sup_{k \geq 1} \Big( \limsup_{n \to \infty} 
\int_Q |T_k (\varrho_n) - T_k(\varrho) |^{\gamma + 1} \ \dxdt \Big) < \infty,
\eF
where $T_k (\varrho) = \min \{ \varrho , k \}$, for any bounded $Q \subset \R^4$ and a $\gamma > 1$.

Due to the poor estimates resulting from \eqref{b31}, \eqref{b32} and \eqref{v8} relation 
\eqref{m9} has to be replaced by a weaker weighted estimate 
\bFormula{m10}
\sup_{k \geq 1} \Big( \limsup_{n \to \infty} 
\int_Q (1 + \vartheta_n )^{- \beta} |T_k (\varrho_n) - T_k(\varrho) |^{\gamma + 1} \ \dxdt \Big) < \infty,
\eF
for suitable $\beta > 0$, $\gamma > 1$.

Finally one recovers strong convergence of the sequence $\{ \vartheta_n \}_{n \geq 1}$ of approximate
temperatures, knowing that spatial gradients of $\vartheta_n$ are uniformly square integrable and that 
\bFormula{m11}
\Ov{ \varrho s(\varrho, \vartheta) \vartheta} = \Ov{ \varrho s (\varrho, \vartheta) } \vartheta,
\eF
where \eqref{m11} can be deduced from the entropy inequality \eqref{v8} by using Aubin-Lions lemma.
\vskip0.5cm
The structure of the proof of Theorem \ref{Tm1} goes as follows: introducing a three level approximation scheme adapted from \cite{EF70} allows to
reduce the proof to a weak stability problem. 
One first proves uniform bounds on the sequence of approximate solutions.
Relying on these uniform bounds, one can handle the convective terms in the field equations by means of a compactness argument.
The strong convergence of the sequence of approximate temperatures is proved by using entropy 
inequality together with the theory of parametrized (Young) measures discussed below.
Then one finally shows that the sequence of approximate densities converges strongly in the Lebesgue space 
$L^1((0,T) \times \Omega)$. This requires weighted estimates of the oscillations defect measure 
in order to show that the limit densities satisfy the renormalized continuity equation. 

\subsubsection{The approximation scheme}

The approximation scheme used in the following is inspired from Chapter 3 in \cite{FEINOV}.

1- The continuity equation \eqref{b9} is regularized by a viscous approximation
\bFormula{a1}
\partial_t \varrho + \Div (\varrho \vc{u}) =  \ep \Delta \varrho  ,\ \ep > 0, 
\eF
 on $(0,T) \times \Omega$, supplemented by  
the homogeneous Neumann boundary conditions
\bFormula{a2}
\Grad \varrho \cdot \vc{n}|_{\partial \Omega} = 0.
\eF
The initial density is given by
\bFormula{a3}
\varrho(0, \cdot) = \varrho_{0, \delta},
\eF
where
\bFormula{a4}
\varrho_{0, \delta} \in C^1(\Ov{\Omega}),\
\Grad \varrho_{0, \delta} \cdot \vc{n} |_{\partial \Omega} = 0,\ 
\inf_{x \in \Omega} \varrho_{0, \delta}(x) > 0, 
\eF  
with a positive parameter $\delta > 0$.

The functions $\varrho_{0,\delta}$ are chosen in such a way that 
\bFormula{a5}
\varrho_{0, \delta} \to \varrho_0 \ \mbox{in}\ L^{5 \over 3}(\Omega),\ 
| \{ \varrho_{0, \delta} < \varrho_0 \} | \to 0 \ \mbox{for}\ 
\delta \to 0.
\eF
Here, of course, the choice of the ``critical'' exponent $\gamma = 5/3$ is directly related 
to estimate \eqref{p7} established above.

2- The regularized momentum equation is
\bFormula{a6}
\partial_t (\varrho \vc{u}) + \Div (\varrho \vc{u} \otimes \vc{u}) + \Grad p +  
\delta \Grad \varrho^{\Gamma} + \ep \Grad \vc{u} \Grad \varrho  = 
\Div (\tn{S}) + \vc{j} \times \vc{B},
\eF
in $(0,T) \times \Omega$ where the approximate velocity field satisfies the no slip boundary condition
\bFormula{a7}
\vc{u}|_{\partial \Omega} = 0. 
\eF
We also prescribe initial conditions
\bFormula{a8}
(\varrho \vc{u})(0, \cdot) = (\varrho \vc{u})_{0, \delta}, 
\eF
where
\bFormula{a9}
(\varrho \vc{u})_{0, \delta} = \left\{ 
\begin{array}{l} (\varrho \vc{u})_0 \ \mbox{provided}\ \varrho_{0, \delta} \geq \varrho_0, \\ \\
0 \ \mbox{otherwise.}
\end{array}
\right.
\eF
The artificial pressure term $\delta \varrho^{\Gamma}$ in \eqref{a6} provides 
additional estimates on the approximate densities in order to facilitate the limit passage 
$\ep \to 0$. To this end, one has to take $\Gamma$ large enough, say, 
$\Gamma > 8$, and to re-parametrize the initial distribution of the approximate densities so that
\bFormula{a10}
\delta \intO{ \varrho_{0, \delta}^{\Gamma} } \to 0 \ \mbox{for}\ \delta \to 0.
\eF

3- The entropy equation \eqref{b16} is replaced by the regularized internal energy balance
\bFormula{a11}
\partial_t (\varrho e + \delta \varrho \vartheta ) + 
\Div \Big( (\varrho e + \delta \varrho \vartheta) \vc{u} \Big) - 
\Div \Big( ( \kappa_F + \kappa_R \vartheta^3 + \delta \vartheta^{\Gamma}) \Grad \vartheta \Big) 
\eF
\[
=\tn{S} : \Grad \vc{u} - p \ \Div (\vc{u}) + \frac{1}{\sigma} \left|\Curl \vc{H}\right|^2 + 
\ep \Gamma |\Grad \varrho|^2 \varrho^{\Gamma - 2},
\]
to be satisfied in $(0,T) \times \Omega$, together with no-flux boundary conditions 
\bFormula{a12}
\Grad \vartheta \cdot \vc{n} |_{\partial} = 0. 
\eF
The initial conditions read
\bFormula{a13}
\varrho (e + \delta \vartheta)(0, \cdot) = \varrho_{0, \delta} (e (\varrho_{0, \delta}, \vartheta_{0, \delta}) + 
\delta \vartheta_{0,\delta} ),
\eF
where the approximate temperature satisfies
\bFormula{a14}
\vartheta_{0, \delta} \in C^{1}(\Ov{\Omega}),\ 
\Grad \vartheta_{0, \delta} \cdot \vc{n}|_{ \partial \Omega} = 0,\ 
\inf_{x \in \Omega} \vartheta_{0, \delta}(x) > 0,
\eF
and
\bFormula{a15}
\left\{ 
\begin{array}{c}
\vartheta_{0, \delta} \to \vartheta_0 \ \mbox{in}\ 
L^p(\Omega) \ \mbox{for any}\ p \geq 1,\\ \\ 
\delta \intO{ \varrho_{0, \delta} \log(\vartheta_{0, \delta}) } \to 0, \\ \\
\intO{ \varrho_{0, \delta} s(\varrho_{0, \delta}, \vartheta_{0, \delta}) } \to 
\intO{ \varrho_0 s(\varrho_0, \vartheta_0 ) }, 
\end{array}
\right. 
\eF
as $\delta \to 0$, with
\bFormula{a16} 
\intO{ \varrho_{0, \delta} e(\varrho_{0, \delta}, \vartheta_{0, \delta}) } < c  
\ \mbox{uniformly for}\ \delta > 0.
\eF
Finally, the magnetic induction vector $\vc{B}$ obeys the original Maxwell's equations
\bFormula{a17}
\partial_t \vc{B} + \Curl (\vc{B} \times \vc{u}) + \Curl \left(\frac{1}{\sigma} \Curl
  \left(\frac{\vc{B}}{\mu}\right) \right) = 0,\ \ \ \Div \vc{B} = 0
\eF
in $(0,T) \times \Omega$, supplemented with the initial condition
\bFormula{a18}
\vc{B}(0, \cdot) = \vc{B}_{0, \delta},
\eF
where, by virtue of Theorem 1.4 in \cite{TEM}, one can take
\bFormula{a19}
\vc{B}_{0, \delta} \in \D(\Omega; \R^3),\ \Div \vc{B}_{0, \delta} = 0,
\eF
\bFormula{a20}
\vc{B}_{0,\delta} \to \vc{B}_0 \ \mbox{in}\ L^2(\Omega; \R^3) \ \mbox{for}\ \delta \to 0.
\eF
For given positive parameters $\ep$, $\delta$ , and $\Gamma > 8$,
the proof of Theorem \ref{Tm1} consists in the following steps:

\begin{description}
\item[Step 1]

Solving problem \eqref{a1}-\eqref{a20} for fixed $\ep > 0$, $\delta > 0$.

\item[Step 2] Passing to the limit for $\ep \to 0$.

\item [Step 3] Letting $\delta \to 0$.

\end{description}

One sees that Step 1 can be achieved by using a simple fixed point argument: equation \eqref{a6} is solved in terms of 
the velocity $\vc{u}$ with the help of the 
Faedo-Galerkin method, where $\varrho$, $\vartheta$, and $\vc{B}$ are computed successively 
from \eqref{a1}, \eqref{a11}, and \eqref{a17} as functions of $\vc{u}$. It is
easy to check that the corresponding approximate solutions satisfy the energy 
balance
\bFormula{a21}
{\rm d \over {\rm d}t} \intO{ 
\Big( {1 \over 2} \varrho |\vc{u}|^2 + \varrho e + {\mathcal M}(|\vc{H}|) + 
 {\delta \over \Gamma - 1} \varrho^{\Gamma} + \delta \varrho \vartheta  \Big) } = 0,
\eF
in $\D'(0,T)$, and the limit property
\bFormula{a22}
\lim_{t \to 0+} \intO{ 
\Big( {1 \over 2} \varrho |\vc{u}|^2 + \varrho e +  {\mathcal M}(|\vc{H}|) + 
 {\delta \over \Gamma - 1} \varrho^{\Gamma} + \delta \varrho \vartheta  \Big) } = 
\eF
\[
\intO{ 
\Big(  {1 \over 2 \varrho_{0, \delta}} | (\varrho \vc{u})_{0, \delta}|^2 + \varrho_{0, \delta} e
(\varrho_{0, \delta}, \vartheta_{0, \delta}) +  {\mathcal M}(|\vc{H}_{0, \delta}|) + 
 {\delta \over \Gamma - 1} \varrho_{0,\delta}^{\Gamma} + 
\delta \varrho_{0, \delta} \vartheta_{0, \delta}  \Big) }.
\]
The reason for splitting  
Steps 2 and 3 into the $\ep$ and $\delta-$parts is a refined density estimates based on multipliers $\Grad (-\Delta)^{-1} [\varrho^{\beta}]$,
where one must take $\beta = 1$ when the artificial viscosity term $\ep \Delta \varrho$ is present, 
while uniform estimates require $\beta$ to be a small positive number. For this reason, we focus only on Step 3:  the goal is
 to establish the  weak sequential stability property for the solution set 
of the approximate problem:

1. The density $\varrho_{\delta} \geq 0$ and the velocity $\vc{u}_{\delta}$ satisfy
\bFormula{a23}
\int_0^T \int_{\Omega}
\Big( \varrho_{\delta} \partial_t \varphi + \varrho_{\delta} \vc{u}_{\delta} \cdot \Grad \varphi \Big) \ \dxdt + 
\int_{\R^3} \varrho_{0, \delta} \varphi (0, \cdot) \ \dx = 0
\eF
for any test function $\varphi \in \D([0,T) \times \R^3)$. In addition, 
\bFormula{a24}
\vc{u}_{\delta}(t) \in W^{1,2}_0 (\Omega; \R^3) \ \mbox{for a.e.}\ t \in (0,T).
\eF
2. The momentum equation 
\bFormula{a25}
\int_0^T \intO{ \Big(
\varrho_{\delta} \vc{u}_{\delta} \cdot \partial_t \varphi + (\varrho_{\delta} \vc{u}_{\delta} \otimes \vc{u}_{\delta}): 
\Grad \varphi + 
p_{\delta} \ \Div (\varphi) +  \delta \varrho^{\Gamma}_{\delta}  \ \Div (\varphi)  \Big) } \ \dt = 
\eF
\[
\int_0^T \intO{ \Big( \tn{S}_{\delta}: \Grad \varphi - 
(\vc{j}_{\delta} \times \vc{B}_{\delta}) 
\cdot \varphi \Big) } \ \dt - \intO{ (\varrho  \vc{u})_{0,\delta} \cdot \varphi (0, \cdot) }
\]
holds for any $\varphi \in \D([0,T) \times \Omega; \R^3)$,
where $p_{\delta} = p(\varrho_{\delta}, \vartheta_{\delta})$, and $\tn{S}_{\delta}$, $\vc{J}_{\delta}$ are determined in terms 
of $\vc{u}_{\delta}$, $\vartheta_{\delta}$, and $\vc{B}_{\delta}$ through the constitutive relations.

3. The entropy production inequality reads now
\bFormula{a27}
\int_0^T \intO{ \Big( (\varrho_{\delta} s_{\delta} +  \delta \varrho_{\delta} \log(\vartheta_{\delta}) ) \partial_t \varphi + 
(\varrho_{\delta} s_{\delta} \vc{u}_{\delta} +  \delta \varrho_{\delta} \log(\vartheta_{\delta}) \vc{u}_{\delta}  ) 
\cdot \Grad \varphi  \Big) } \ \dt
\eF
\[
+\int_0^T \intO{\left( {\vc{q}_{\delta} \over \vartheta_{\delta}} \cdot \Grad \varphi_{\delta} -  \delta \vartheta^{\Gamma - 1}_{\delta} 
\Grad \vartheta_{\delta}  \cdot \Grad \varphi \right)} \ \dt  
\]
\[
\leq \int_0^T \intO{ {1 \over \vartheta_{\delta}} \Big( {\vc{q}_{\delta} \cdot \Grad \vartheta_{\delta} \over \vartheta_{\delta} }
 - \tn{S}_{\delta} : \Grad \vc{u}_{\delta} - {1\over\sigma_{\delta}} \left|\Curl \vc{H}_{\delta}\right|^2
 -  \delta \vartheta^{\Gamma - 1}_{\delta} |\Grad \vartheta_{\delta}|^2   \Big) \varphi  } \ \dt 
\]
\[
-\intO{ \Big( \varrho_{0, \delta} s(\varrho_{0, \delta}, \vartheta_{0, \delta})
+  \delta \varrho_{0,\delta} \log(\vartheta_{0,\delta})  \Big) \varphi (0, \cdot) },
\]
for any $\varphi \in \D([0,T) \times \R^3)$, $\varphi \geq 0$. Here, $\vartheta_{\delta}$ is assumed to be positive a.e. on the set
$(0,T) \times \Omega$, $s_{\delta} = s(\varrho_{\delta}, \vartheta_{\delta})$, $\sigma_{\delta} = 
\sigma (\varrho_{\delta}, \vartheta_{\delta})$, and $\vc{q}_{\delta}$ is a function of $\varrho_{\delta}$, 
$\vartheta_{\delta}$ given by \eqref{b28}.

4. The magnetic induction vector $\vc{B}_{\delta}$ satisfies
\bFormula{a28}
\int_0^T \intO{ \Big( \vc{B}_{\delta} \cdot \partial_t \varphi - ( \vc{B}_{\delta} \times \vc{u}_{\delta} + 
{1\over\sigma_{\delta}} \ \Curl \vc{H}_{\delta} ) \cdot \Curl \varphi \Big) } \ \dt + 
\intOd{ \vc{B}_{0, \delta} \cdot \varphi (0, \cdot) } = 0, 
\eF
for any vector field $\varphi \in \D([0,T) \times \R^3; \R^3)$. 

5. Finally, the (total) energy equality
\bFormula{a29}
\intO{ \Big( {1 \over 2} \varrho_{\delta} |\vc{u}_{\delta}|^2 + \varrho_{\delta} e + {\mathcal M}(\vc{H}_{\delta}) 
 +  {\delta \over \Gamma - 1} 
\varrho^{\Gamma}_{\delta} + \delta \varrho_{\delta} \vartheta_{\delta}  \Big)(t) }  
\eF
\[
=\intO{ 
\left( {1 \over 2 \varrho_{0, \delta}} | (\varrho \vc{u})_{0, \delta}|^2 + \varrho_{0, \delta} e
(\varrho_{0, \delta}, \vartheta_{0, \delta}) + {\mathcal M}(\vc{H}_{0,\delta})\right) } +
\intOd{  \left(
{\delta \over \Gamma - 1} \varrho_{0,\delta}^{\Gamma} + 
\delta \varrho_{0, \delta} \vartheta_{0, \delta}  \right)},
\]
holds for a.e. $t \in (0,T)$.

The weak sequential stability problem to be addressed below consists now in showing 
that one can pass to the limit 
\[
\varrho_{\delta} \to \varrho,\ \vc{u}_{\delta} \to \vc{u},\ \vartheta_{\delta} \to \vartheta, \ \vc{B}_{\delta} \to \vc{B},
\]
as $\delta \to 0$,
in a suitable topology, where the limit quantity $ \{ \varrho, \vc{u}, \vartheta, \vc{B} \}$ is a variational 
solution of problem \eqref{v1}-\ref{v15}.

\subsubsection{Uniform bounds}
\label{u}

One must identify uniform bounds on the sequences $\sD{ \rD }$, $\sD{ \uD }$, 
$\sD{ \tD }$, and $\sD{ \bD }$ through total energy balance \eqref{a29}, dissipation inequality 
\eqref{a27} and the other relations resulting from \eqref{a23}-\eqref{a29}. 
\vskip0.25cm
{\it 1. Total mass conservation}. As $\rD$, $\rD \uD$ are locally integrable in $[0,T) \times \Ov{\Omega}$, it follows from 
\eqref{a23} that total mass is a constant of motion, specifically, 
\bFormula{u1}
\intO{ \rD (t) } = \intO{ \varrho_{0, \delta} } = M_0 \ \mbox{for a.e.}\ t \in (0,T).
\eF
In particular, as $\rD \geq 0$ and \eqref{a5} holds, we get 
\bFormula{u2}
\sD{ \rD } \ \mbox{bounded in}\ L^{\infty}(0,T; L^1(\Omega)).
\eF
{\it 2. Energy estimates}. Using \eqref{b25} and \eqref{p7} one gets
\[
\intO{ \rD e (\rD, \tD) } = {3 \over 2} \intO{ p_F(\rD , \tD) } \geq {3 p_{\infty} \over 2} 
\intO{ \rD^{ 5 \over 3} }.
\]  
Total energy balance \eqref{a29} together with bounds on the initial data 
\eqref{a5}, \eqref{a9}, \eqref{a10}, \eqref{a16}, and \eqref{a20} yield the energy estimates
\bFormula{u4}
\sD{ \rD } \ \mbox{bounded in}\ L^{\infty}(0,T; L^{5 \over 3}(\Omega)),
\eF
\bFormula{u5}
\sD{ \tD } \ \mbox{bounded in}\ L^{\infty}(0,T; L^4(\Omega)),
\eF
\bFormula{u6}
\sD{ \vc{H}_{\delta} } \ \mbox{bounded in}\ L^{\infty}(0,T; L^2(\Omega; \R^3)), 
\eF
\bFormula{u7}
\sD{ \rD |\uD |^2 } ,\ 
\sD{ \rD e(\rD, \tD ) } \ \mbox{bounded in}\ L^{\infty}(0,T; L^1(\Omega)),
\eF
and
\bFormula{u8}
\delta \intO{ \rD^{\Gamma} } \leq c \ \mbox{uniformly with respect to}\ \delta > 0.
\eF
In particular, using H\"older's inequality, \eqref{u4} and \eqref{u7} imply
\bFormula{u9}
\sD{ (\varrho \vc{u})_{\delta} } \ \mbox{bounded in}\ L^{\infty}(0,T; L^{5 \over 4}(\Omega; \R^3)).
\eF
{\it 3. Dissipation estimates}. Taking a spatially homogeneous test function $\varphi$ such that $\varphi(0, \cdot) = 1$ in the entropy 
production inequality \eqref{a27} and using \eqref{a15} \eqref{a16},
we obtain
\bFormula{u10}
\int_0^\tau \intO{ {1 \over \tD} \Big( \tn{S}_{\delta} : \Grad \uD  - {\vc{q}_{\delta} \cdot \Grad \tD \over \tD } 
 + {1\over \sigma_{\delta}} |\Curl \vc{H}_{\delta} |^2 + 
\delta \tD^{\Gamma - 1} |\Grad \tD |^2 \Big)} \ \dt 
\eF
\[
\leq c_1 +  \intO{ \rD(\tau) s(\rD(\tau), \tD(\tau)) + \delta \rD(\tau) \log( \tD(\tau)) } 
\leq 
c_2 + \intO{ \rD(\tau) s(\rD(\tau), \tD(\tau)) },
\]
for a.e.$\ \tau \in (0,T)$, where the last inequality follows from \eqref{u4}, \eqref{u5}. 

Using \eqref{p12} and \eqref{p13}, the integral in the right hand side may be estimated
\begin{multline}
  \label{u11}
\intO{ \rD s (\rD, \tD) } = \intO{ {4 \over 3} a \tD^3 + \rD S_F \left( {\rD \over \tD^{3 \over 2} } \right) } 
\leq c_1 - c_2 \int_{ \{ \rD \leq \tD^{3 \over 2} \} } \rD \left( \log(\rD) - {3 \over 2} \log(\tD) \right) \ \dx
 \\ \leq 
c_3 + c_4 \int_{ \{ \rD \leq \tD^{3 \over 2} \} } \rD \log( \tD ) \ \dx \leq c_5 + c_6 \intO{ \rD \tD } \leq c_7,  
\end{multline}
where we have taken into account \eqref{u4} and \eqref{u5}.

Thus the integral in the left hand side of \eqref{u10} is bounded independently of $\delta$ and we conclude after
\eqref{b29}-\eqref{b32bis} that
\bFormula{u12}
\sD{ (1 + \tD )^{\alpha - 1 \over 2} \left\langle \Grad \uD \right\rangle }\ \mbox{is bounded in}\ 
L^2(0,T; L^2(\Omega; \R^{3 \times 3})), 
\eF
\bFormula{u12bis}
\sD{ (1 + \tD )^{\alpha - 1 \over 2} \Div \left(\uD\right) } \ \mbox{is bounded in}\ 
L^2(0,T; L^2(\Omega)),
\eF
\bFormula{u13}
\sD{ \Grad \log( \tD) },\  
\sD{ \Grad \tD^{3 \over 2} } \ \mbox{are bounded in}\ L^2(0,T; L^2(\Omega; \R^3)),
\eF
\bFormula{u14}
\sD{ (1 + \tD)^{\alpha - 1 \over 2} \Curl \vc{H}_{\delta} } \ \mbox{is bounded in}\ L^2(0,T; L^2(\Omega; \R^3)),
\eF
and
\bFormula{u15}
\sD{ \sqrt{ \delta } \Grad \tD^{\Gamma \over 2} } \ \mbox{is bounded in}\ L^2(0,T; L^2(\Omega; \R^3)),
\eF
where $\alpha$ satisfies \eqref{m1}, and where we denote by
\[
\left\langle \tn{Q} \right\rangle = {1 \over 2}(\tn{Q} + \tn{Q}^T ) - {1 \over 3} {\rm trace}(\tn{Q}) \tn{I}, 
\]
the traceless component of the symmetric part of a tensor $\tn{Q}$.

Since the velocity field satisfies \eqref{a24} we find 
\bFormula{u16}
\sD{ \uD } \ \mbox{bounded in}\ L^2(0,T; W^{1,2}_0 (\Omega; \R^3)),
\eF
and \eqref{u13} and \eqref{u5} lead to 
\bFormula{u17}
\sD{ \tD^{3 \over 2} } \ \mbox{bounded in}\ L^2(0,T; W^{1,2}(\Omega)). 
\eF
{\it 4. Positivity of the absolute temperature}. According to physics, temperature must be positive a.e. on $(0,T) \times \Omega$.
 To justify that we use the 
uniform $L^2-$ estimates of $\Grad \log(\tD)$ given by \eqref{u16} and the following version
of Poincar\'e's inequality:
\bLemma{u1}
Let $\Omega \subset \R^N$ be a bounded Lipschitz domain, and $\omega \geq 1$ be a given constant. Furthermore, assume 
$O \subset \Omega$ is a measurable set such that $| O | \geq \underline{o} > 0$. Then 
\[
\| v \|_{W^{1,2}(\Omega)} \leq c(\underline{o}, \Omega, \omega) \Big( \| \Grad v \|_{L^2(\Omega; \R^3)} + 
( \int_O |v|^{1 \over \omega} \ \dx )^{\omega} \Big).
\]
\eL
In order to apply Lemma \ref{Lu1}, we show first the uniform bound
\bFormula{u18}
0 < \underline{T} \leq \intO{ \tD (t) } \ \mbox{for a.e.}\ t \in (0,T).
\eF
As a consequence of \eqref{a27} we get
\[
\intO{ \Big(\rD s (\rD, \tD) + \delta \rD \log( \tD )\Big)(t) } \geq \intO{ \Big(
\varrho_{0, \delta} s (\varrho_{0, \delta}, \vartheta_{0, \delta}) + \delta \varrho_{0, \delta} \log(\vartheta_{0, \delta})\Big) },
\]
where, by virtue of \eqref{a5} and \eqref{a15} 
\[
\delta \intO{ \varrho_{0, \delta} \log( \vartheta_{0, \delta}) } \to 0\ \ \mbox{and}\ \ 
\intO{ \varrho_{0, \delta} \ s(\varrho_{0, \delta}, \vartheta_{0, \delta}) } \to 
\intO{ \varrho_0 s (\varrho_0, \vartheta_0 ) }.
\]
On the other hand \eqref{u4} and \eqref{u5} yield
\[
\delta \intO{ \rD \log (\tD )(t) } \leq 
\delta \intO{ \rD(t) \tD(t) } \to 0\ \ \ \ \mbox{for a.e}\ t \in (0,T),
\]
while, using \eqref{p9}
\[
\intO{ \rD s( \rD , \tD )(t) } = 
\intO{ \left({4 \over 3} a \tD^3  + \rD S_F \left( { \rD \over \tD^{3 \over 2} } \right)\right)(t) }.
\]
Supposing the opposite of \eqref{u18}, one could extract a sequence 
$ \sD{ \tD (t_{\delta}) }$ such that
\bFormula{u19}
\tD (t_{\delta}) \to 0 \ \mbox{weakly in}\ L^4(\Omega), 
\ \mbox{and strongly in}\ L^p(\Omega) \ \mbox{for any}\ 1 \leq p < 4,
\eF
\bFormula{u20}
\intO{\left( {4 \over 3} a \vartheta_0^3 + \varrho_0 S_F\left( { \varrho_0 \over \vartheta_0^{3 \over 2}} \right) \right)}
\leq \liminf_{ \delta \to 0+ } \intO{ \rD (t_{\delta}) S_F \Big( {\rD \over \tD^{3 \over 2}} \Big) (t_{\delta}) },
\eF 
where
\[
\rD (t_{\delta}) \to \varrho (t) \ \mbox{weakly in}\ L^{5 \over 3}(\Omega),\ 
\intO{ \rD (t_{\delta}) } = \intO{ \varrho_{0, \delta} }.
\]
Now, for any fixed $K > 1$, one can write
\bFormula{u21}
\intO{ \rD (t_{\delta}) S_F \Big( { \rD \over \tD^{3 \over 2} } \Big) (t_{\delta}) }
\eF
\[
\leq\int_{ \{ \rD(t_{\delta}) \leq K \tD^{3 \over 2}(t_{\delta})  \} } \rD (t_{\delta}) S_F \left( 
{ \rD \over \tD^{3 \over 2} } \right) (t_{\delta})\ \dx + 
\int_{ \{ \rD(t_{\delta}) \geq K \tD^{3 \over 2}(t_{\delta})  \} } \rD (t_{\delta}) S_F \left( 
{ \rD \over  \tD^{3 \over 2} } \right) (t_{\delta}) \ \dx 
\]
\[\leq
c \int_{\{ \rD(t_{\delta}) \leq K \tD^{3 \over 2}(t_{\delta})  \} } \rD (t_{\delta}) \left( 1 + {\tD^{3 \over 2} \over \rD } (t_{\delta}) \right) \ \dx
+ S_F(K) \intO{ \varrho_{0, \delta} }
\]
\[
 \leq 2 c \intO{ \tD^{3 \over 2}( t_{\delta} ) } + S_F(K) \intO{ \varrho_{0, \delta} }. 
\]
Thus combining \eqref{u19}-\eqref{u21} one concludes
\[
\intO{ {4 \over 3} a \vartheta_0^3 + \varrho_0 S_F\Big( { \varrho_0 \over \vartheta_0^{3 \over 2}} \Big) }
\leq S_F (K) \intO{ \varrho_0 } \ \mbox{for any}\ K > 1,
\]
which is in contradiction with the fact that $S_F$ is strictly decreasing with $\lim_{K \to \infty} S_F(K) < 0$, and $\vartheta_0$ positive (non-zero) on 
$\Omega$. This proves the uniform bound \eqref{u18}. 

Finally, observing that 
\[
\underline{T} - \ep |\Omega| \leq \int_{ \{ \tD (t) > \ep \} } \tD(t) \ \dx \leq 
| \{ \tD(t) > \ep \} |^{3 \over 4} \| \tD(t) \|_{L^4(\Omega)} \ \mbox{for any}\ \ep > 0,
\]
we deduce, using \eqref{u5}, that there exist $\ep > 0$ and $\underline{o}$ such that
\[
| \{ \tD(t) > \ep \} | > \underline{o} > 0 \ \mbox{for a.e.}\ t \in (0,T),
\]
uniformly for $\delta > 0$. Then using \eqref{u13} one can apply Lemma \ref{Lu1} to $\log (\tD)$ and we obtain 
the required estimate
\bFormula{u22}
\sD{ \log (\tD ) } \ \mbox{bounded in}\ L^2(0,T; W^{1,2}(\Omega)).
\eF 
{\it 5. Estimates of the magnetic field}. Using the integral identity \eqref{a28} we get
\bFormula{u23}
\Div (\bD) (t) = 0 \ \mbox{in}\ \D'(\Omega),\ \bD \cdot \vc{n}|_{\partial \Omega} = 0,
\eF
which, together with estimates \eqref{u6}, \eqref{u14} and Theorem 6.1 in \cite{DULI}, yields
\bFormula{u24}
\sD{ \vc{H}_{\delta} } \ \mbox{bounded in}\ L^2(0,T; W^{1,2}(\Omega; \R^3)).
\eF
{\it 6. Pressure estimates}. We first observe that estimates \eqref{u4}, \eqref{u16} imply boundedness of 
the sequences $\sD{ \rD \uD }$ and $\sD{ \rD \uD \otimes \uD }$ in 
$L^p((0,T) \times \Omega)$ for a certain $p > 1$. 

Now one gets
\bFormula{u25}
\tn{S}_{\delta} = \sqrt{ \tD \nu (\tD) } \sqrt{ {\nu(\tD) \over \tD}} \left\langle \Grad \uD \right\rangle + 
\sqrt{ \tD \eta (\tD)} \sqrt{ {\eta(\tD) \over \tD }} \Div( \uD), 
\eF
where, by virtue of \eqref{u5}, \eqref{u12}, \eqref{u17}, and hypothesis \eqref{m1}, the expression on 
the right-hand side is bounded in $L^p((0,T) \times \Omega)$ for a certain $p > 1$. 

Finally the Lorentz force $\vc{j}_{\delta} \times \bD=\Curl\vc{H}_{\delta} \times \mu(\vc{H}_{\delta})\vc{H}_{\delta}$ is controlled by 
\eqref{u6}, \eqref{u24} and \eqref{b32bis} which, combined with a simple interpolation argument, yields
\bFormula{u26}
\sD{ \vc{j}_{\delta} \times \bD } \ \mbox{bounded in}\ L^p((0,T) \times \Omega)
\ \mbox{for a certain}\ p > 1.
\eF
At this stage, using the main result of \cite{FP13}, we are allowed to use the quantities
\[
\varphi(t,x) = \psi(t) {\cal B}[ \rD^{\omega} ],\ \psi \in \D(0,T)\ \ \mbox{for a sufficiently small parameter}\ 
\omega > 0,
\]
as test functions in \eqref{a25}, where ${\cal B}$ is the Bogovskii operator \cite{BOG} solving the boundary value problem 
\bFormula{u27}
\Div \Big( {\cal B} [v] \Big) = v - {1 \over |\Omega|} \intO{ v },\ \ \ {\cal B}|_{\partial \Omega} = 0.
\eF
The resulting estimate reads
\bFormula{u28}
\int_0^T \intO{ \Big( p(\rD , \tD) \rD^{\omega} + \delta \rD^{\Gamma + \omega} \Big) } \ \dt < c, 
\ \mbox{with}\ c \ \mbox{independent of}\ \delta,
\eF
then
\bFormula{u29}
\sD{ p (\rD, \tD) } \ \mbox{is bounded in}\ L^p ((0,T) \times \Omega)
\ \mbox{for a certain}\ p > 1,
\eF
and
\bFormula{u30}
\sD{ \rD^{{5 \over 3} + \omega} } \ \mbox{is bounded in}\ 
L^1((0,T) \times \Omega).
\eF
\subsubsection{Sequential stability of the equations}

 1. {\it Continuity equation}. With the estimates established in the previous section, we can pass to the limit 
 $\delta \to 0$ in \eqref{a23}.

 Indeed, extracting subsequences if necessary, we deduce from \eqref{u4},\eqref{u9}, 
\eqref{u16}, and the fact that $\rD$ and $\uD$ satisfy the identity \eqref{a23}
\bFormula{s1}
\rD \to \varrho \ \mbox{in}\ C([0,T]; L^{5 \over 3}_{weak} (\Omega)),
\eF
\bFormula{s2}
\uD \to \vc{u} \ \mbox{weakly in}\ L^2(0,T; W^{1,2}_0 (\Omega; \R^3)),
\eF
\bFormula{s3}
\rD \uD \to \varrho \vc{u} \ \mbox{weakly-(*) in}\ L^{\infty}(0,T; L^{5 \over 4}(\Omega; \R^3)),
\eF
where the limit quantities satisfy \eqref{v1}. 
\vskip0.25cm
 2.{\it Momentum equation}. Using the estimates obtained in Section \ref{u}.7 together with \eqref{a25} we have
\bFormula{s4}
\rD \uD \to \varrho \vc{u} \ \mbox{in}\ C([0,T]; L^{5 \over 4}(\Omega; \R^3)),
\eF
\bFormula{s5}
\rD \uD \otimes \uD \to 
\varrho \vc{u} \otimes \vc{u} \ \mbox{weakly in}\ L^2(0,T; L^{30 \over 29}(\Omega; \R^{3 \times 3}_{sym})),
\eF
where we have used the embedding $W^{1,2}_0 (\Omega) \hookrightarrow L^6(\Omega)$.

Furthermore, in accordance with \eqref{u28}, \eqref{u29}
\bFormula{s6}
p(\rD, \tD) \to \Ov{ p(\varrho, \vartheta) } \ \mbox{weakly in}\ 
L^p((0,T) \times \Omega), 
\eF
and
\bFormula{s7}
\delta \rD^{\Gamma} \to 0 \ \mbox{in}\ L^p((0,T) \times \Omega),
\eF
for a certain $p > 1$. Here, in agreement with Section \ref{u}, 
\[
\int_0^T \intO{ \Ov{ p(\varrho, \vartheta) } \varphi } \ \dt = 
\int_0^T \intO{ \varphi \Big( \int_{\R^2} p(\varrho, \vartheta) \ {\rm d}\Lambda_{t,x} (\varrho, \vartheta) \Big) }
\ \dt ,\ \varphi \in \D((0,T) \times \Omega),
\]
where $\Lambda_{t,x}(\varrho, \vartheta)$ is a parametrized (Young) measure associated to the (vector valued) 
sequence $\sD{ \rD , \tD }$.

Similarly, one can use \eqref{a28} together with estimates \eqref{u6}, \eqref{u24} to deduce
\bFormula{s8}
\bH \to \vc{H} \ \mbox{weakly in}\ 
L^2(0,T; W^{1,2}(\Omega; \R^3)) \ \mbox{and strongly in}\ L^2((0,T) \times \Omega; \R^3),
\eF
and
\[
\mu(\bH) \to \mu(\vc{H}) \  \ \mbox{strongly in}\ L^2((0,T) \times \Omega; \R^3),
\]
 then
\[
\vc{j}_{\delta} \times \vc{B}_{\delta}\ \to \vc{j} \times \vc{B}\ \mbox{weakly in}\ L^p((0,T) \times \Omega; \R^3),
\]
for a certain $p > 1$. Thus the limit quantities satisfy an averaged momentum equation 
\bFormula{s9}
\int_0^T \intO{ 
\Big( (\varrho \vc{u}) \cdot \partial_t \varphi + (\varrho \vc{u} \otimes \vc{u}): \Grad \varphi + 
\Ov{p (\varrho , \vartheta )} \ \Div (\varphi) \Big) } \ \dt 
\eF
\[
=\int_0^T \intO{ \Big( 
\Ov{ \tn{S} }: \Grad \varphi - (\vc{j} \times \vc{B} ) \cdot \varphi \Big)
} \ \dt - \intO{ (\varrho \vc{u})_0 \cdot \varphi (0, \cdot) }
\]
for any vector field $\varphi \in \D([0,T) \times \Omega; \R^3)$.

The symbol $\Ov{S}$ denotes a weak limit in $L^p(0,T; L^p( \Omega; \R^{3 \times 3}_{sym}))$, $p > 1$ of the approximate 
viscosity tensors $\tn{S}_{\delta}$ specified in \eqref{u25}. 
Clearly, relation \eqref{s9} will coincide with the (variational) momentum 
equation \eqref{v4} as soon as we show strong (pointwise) convergence of the sequences 
$\sD{ \rD }$ and $\sD{ \tD }$.
\subsubsection{ Entropy inequality and strong convergence of the temperature}

1. {\it Entropy inequality}. In order to obtain information on the time oscillations of the sequence 
$\sD{ \tD }$, we use estimates on $\partial_t (\rD s(\rD, \tD))$
provided by the approximate entropy balance \eqref{a27}.

It first follows from \eqref{p12}, \eqref{p13} that
\[
|\rD s(\rD, \tD)| \leq c \Big( \tD^3 + \rD | \log (\rD ) | + \rD |\log (\tD )| \Big),
\]
therefore, by virtue of the uniform estimates \eqref{u4}, \eqref{u5}, \eqref{u9}, and \eqref{u22}, we can assume
\bFormula{e1}
\rD s (\rD, \tD ) \to 
\Ov{ \varrho s(\varrho, \vartheta) } \ \ \ \mbox{weakly in} \ L^p((0,T) \times \Omega),
\eF
\bFormula{e2}
\rD s( \rD, \tD) \uD \to \Ov{ \varrho s(\varrho, \vartheta) \vc{u} }
\ \ \ \mbox{weakly in}\ L^p((0,T) \times \Omega; \R^3),
\eF
for a certain $p > 1$.

Similarly, one estimates the entropy flux
\bFormula{e3}
\Big| {\vc{q} \over \tD } + \delta \tD^{\Gamma - 1} \Grad \tD \Big| \leq c 
\Big( |\Grad \log(\tD) | + \tD^2 |\Grad \tD | + \delta \tD^{\Gamma - 1} |\Grad \tD | \Big),
\eF
where 
\[
\tD^2 |\Grad \tD| = {2 \over 3} \tD^{3 \over 2} \left|\Grad \tD^{3 \over 2}\right|,\ 
 \delta \tD^{\Gamma - 1} \Grad \tD = {2 \over \Gamma} \sqrt{\delta} \ \left| \Grad \tD^{\Gamma \over 2} \right| 
\ \sqrt{\delta} \tD^{s {\Gamma \over 2}} \tD^{(1 - s){\Gamma \over 2}}. 
\]
Choosing the parameter $0 < s < 1$ small enough so that
\[
(1 - s) {\Gamma \over 2} = 1,
\]
we have after H\"older's inequality  
\[
\left\| \ \delta \tD^{\Gamma - 1} \Grad \tD \ \right\|_{L^p(\Omega; \R^3)} \leq 
c \left\| \ \sqrt{\delta} \Grad \tD^{\Gamma \over 2} \ \right\|_{L^2(\Omega; \R^3)} \left\| \tD \right\|_{L^4(\Omega)}
\left\| \ \sqrt{\delta} \tD^{s {\Gamma \over 2}} \ \right\|_{L^6(\Omega)}.
\]
Then using estimates \eqref{u5}, \eqref{u15} together with
the imbedding $W^{1,2}(\Omega) \hookrightarrow L^6(\Omega)$, we get
\bFormula{e4}
\delta \tD^{ \Gamma - 1} \Grad \tD \to 0 
\ \mbox{in}\ L^p((0,T) \times \Omega; \R^3) \ \mbox{for a certain}\ p > 1.
\eF
Moreover using similar arguments one can also show that 
\bFormula{e5}
\sD{ {\vc{q} \over \tD } } \ \mbox{is bounded in}\ 
L^p((0,T) \times \Omega; \R^3),\ \mbox{for a certain}\ p > 1.
\eF
 On the other hand according to \eqref{u13}
\bFormula{e6}
b(\tD) \to \Ov{ b(\vartheta) } \ \mbox{weakly in}\ 
L^q ((0,T) \times \Omega)),\ \mbox{and weakly in}\ 
L^2(0,T; W^{1,2}(\Omega)), 
\eF
for any finite $q \geq 1$ and any function $b$ provided both $b$ and $b'$ are uniformly bounded.

Now, as a consequence of the entropy balance \eqref{a27}, we have 
\[
{\rm Div}_{t,x} \Big[ \rD s(\rD, \tD) + \delta \rD \log( \tD),\ 
( \rD s(\rD, \tD) + \delta \rD \log( \tD)) \uD + {\vc{q} \over \tD} - \delta \tD^{\Gamma - 1} 
\Grad \tD \Big] 
\]
\[
\geq 0 \ \ \ \ \mbox{in} \ \D'((0,T) \times \Omega),
\]
while \eqref{e6} yields
\[
{\rm Curl}_{t,x} \Big[ b(\tD ), 0 , 0 , 0 \Big] \ \mbox{bounded in}\ L^2((0,T) \times \Omega; \R^{4 \times 4}).
\]
Thus a direct application of the Div-Curl lemma \cite{M,T} gives the relation
\bFormula{e7}
\Ov{ \varrho s (\varrho, \vartheta ) b (\vartheta) } = \Ov{ \varrho s (\varrho, \vartheta) } \ \Ov{b (\vartheta)},
\eF
for any $b$ as in \eqref{e6}. 

Moreover, seeing that the sequence $\sD{ \rD s(\rD, \tD) \tD}$ is bounded in 
$L^p((0,T) \times \Omega)$ we deduce from \eqref{e8} the formula
\bFormula{e8} 
\Ov{ \varrho s (\varrho, \vartheta ) \vartheta } = \Ov{ \varrho s (\varrho, \vartheta) } \vartheta.
\eF
2. {\it Parametrized measures and pointwise convergence of the temperature}
\vskip0.25cm
To show that \eqref{e8} implies strong convergence of the sequence  $\sD{ \tD }$, we first remark that the approximate 
solutions solve the renormalized continuity equation 
\bFormula{e9}
\partial_t b(\rD) + \Div ( b(\rD) \uD ) + \Big( b'(\rD) \rD - b(\rD) \Big) \ \Div (\uD) = 0
\ \ \ \mbox{in}\ \D'((0,T) \times \R^3),
\eF
provided that $\rD$ and $\uD$ are extended by zero outside $\Omega$, and for any continuously differentiable function $b$ whose 
derivative vanishes for large arguments.

Functions $\rD$ being square-integrable because of
the artificial pressure in the energy equality, \eqref{e9} follows from \eqref{a23} via the technique of 
DiPerna and Lions \cite{DL}.

Now it follows from \eqref{e9} that 
\bFormula{e10}
b(\rD) \to \Ov{ b(\varrho) } \ \mbox{in}\ C([0,T]; L^p_{weak}(\Omega))\ \ \ \mbox{for any finite}\ p > 1 
\ \mbox{and bounded}\ b.
\eF
Relations \eqref{e10} and \eqref{e6} yield
\[
\Ov{ g(\varrho) h (\vartheta) } = \Ov{ g(\varrho)} \ \Ov{ h(\vartheta) },
\]
or, in terms of the corresponding parametrized measures (see Section \ref{param})
\bFormula{e11}
\Lambda_{t,x}(\varrho, \vartheta) = \Lambda_{t,x}(\varrho) \otimes \Lambda_{t,x} (\vartheta).
\eF
Relation \eqref{e11} means that oscillations in the sequences $\sD{ \rD }$ and 
$\sD{ \tD }$ are orthogonal, in the sense that the parametrized measure associated to $\sD{ \rD, \tD }$ can be written as a tensor product
of the parametrized measures generated by $\sD{ \rD }$ and $\sD{ \tD }$. 

Let us consider a function $H = H(t,x, r, z)$ defined for $t \in (0,T)$, $x \in \Omega$, and $(r,z) \in \R^2$ through formula
\[
H(t,x,r,z) = r \Big( s (r,z) - s(r, \vartheta(t,x) \Big) \Big( z - \vartheta(t,x) \Big).
\]

Clearly, $H$ is a Caratheodory function: $H(t,x, \cdot, \cdot)$ is continuous for a.e $(t,x) \in (0,T) \times 
\Omega$, and $H(\cdot, \cdot, r, z)$ is measurable for any $(r,z) \in \R^2$.

Moreover, as both $s_R$ and $s_F$ are increasing functions of the absolute temperature, we have
\bFormula{e12}
H(t,x,r,z) \geq {4 \over 3} a \Big( z^3 - \vartheta^3(t,x) \Big) \Big( z - \vartheta (t,x) \Big) \geq 0.
\eF
At this stage, we use Theorem 6.2 in \cite{PED1}, namely: weak limits of Caratheodory functions can be expressed in 
terms of the parametrized measure. Accordingly, we obtain
\[
\lim_{\delta \to 0} 
\int_0^T \intO{ \varphi(t,x) H(t,x, \rD, \tD ) } \ \dt = 
\int_0^T \intO{ \varphi(t,x) \Big( \int_{\R^2} H(t,x, \varrho, \vartheta) \ {\rm d} \Lambda_{t,x}(\varrho, \vartheta) \Big) } \ \dt 
\]
\[
=\int_0^T \intO{ \varphi(t,x) \Big( \int_{\R^2} \varrho s (\varrho, \vartheta)(\vartheta - \vartheta(t,x)) {\rm d}
 \Lambda_{t,x}(\varrho, \vartheta) \Big) } \ \dt  
\]
\[
-\int_0^T \intO{ \varphi(t,x) \Big( \int_{\R^2} \varrho s (\varrho, \vartheta(t,x))(\vartheta - \vartheta(t,x)) {\rm d}
 \Lambda_{t,x}(\varrho, \vartheta) \Big) } \ \dt 
\]
\[
=\int_0^T \intO{ \varphi(t,x) \Big( \int_{\R^2} \varrho s (\varrho, \vartheta)(\vartheta - \vartheta(t,x)) {\rm d}
 \Lambda_{t,x}(\varrho, \vartheta) \Big) } \ \dt 
\]
\[
=\int_0^T \intO{ \varphi \Big( \Ov{ \varrho s (\varrho , \vartheta) \vartheta } - 
\Ov{\varrho s (\varrho, \vartheta) } \vartheta \Big) } \ \dt = 0 \ \mbox{for any}\ \varphi \in \D((0,T) \times \Omega,
\] 
where we have used \eqref{e8} to get the last equality together with \eqref{e11}, observing that 
\[
\int_{\R^2} \varrho s (\varrho, \vartheta(t,x)) ( \vartheta - \vartheta(t,x)) \ {\rm d}\Lambda_{t,x}(\varrho, \vartheta) = 
\int_R \varrho s (\varrho, \vartheta (t,x)) \ {\rm d}\Lambda_{t,x}(\varrho) 
\int_R (\vartheta - \vartheta(t,x)) \ {\rm d}\Lambda_{t,x}(\vartheta) = 0.
\]
In particular, we deduce from \eqref{e12} that  $\Ov{ \vartheta^3  \vartheta } = \Ov{\vartheta^3 } \vartheta$, which is equivalent to 
\bFormula{e13}
\tD \to \vartheta \ \mbox{in}\ L^4((0,T) \times \Omega).
\eF

\subsubsection{Pointwise convergence of densities}
\label{d}

{\it 1. The effective viscous pressure}. Let us start with a result of P.-L. Lions \cite{LI4} on the effective viscous 
pressure.

 In the present setting, it can be stated in terms of the parametrized measures as follows
\bFormula{d1}
\Ov{ \psi p (\varrho, \vartheta) b(\varrho)} - \Ov{ \psi p(\varrho, \vartheta)}\ \Ov{b(\varrho)} =
\Ov{ {\cal R}: [\psi \ \tn{S}] b(\varrho) } - \Ov{ {\cal R}: [\psi \ \tn{S}] } \ \Ov{b(\varrho)},
\eF 
for any $\psi \in \D(\Omega)$, and any bounded continuous function $b$, where
${\cal R} = {\cal R}_{i,j}$ is the Riesz operator defined by means of the Fourier transform:
\bFormula{d2}
{\cal R}_{i,j}[v] = \partial_{x_i} \Delta^{-1}_x \partial_{x_j} v  = 
{\cal F}^{-1}_{\xi \to x} \Big[ {\xi_i \xi_j \over |\xi |^2 } {\cal F}_{x \to \xi}[v] \Big].
\eF
Note that \eqref{d1} does not depend on the specific form of the constitutive relations for $p$, $\tn{S}$ and requires only 
validity of the momentum equation \eqref{a25} and the renormalized continuity equation \eqref{e9} for 
$\rD$ and $\uD$.
\vskip0,25cm
{\it 2. Commutator estimates}. If viscosity coefficients $\nu$ and $\eta$ were constant, we would have ${\cal R}[ \tn{S} ] = ({4 \over 3} \nu + \eta ) \ \Div (\vc{u})$, 
and \eqref{d1} would give immediately
\bFormula{d3}
\Ov{p (\varrho, \vartheta) b(\varrho)} - \Ov{p(\varrho, \vartheta)}\ \Ov{b(\varrho)} =
\left({4 \over 3} \nu + \eta\right) \Ov{ \Div (\vc{u} )\  b(\varrho) } - \left({4 \over 3} \nu + \eta\right) \Div(
\vc{u})  \ \Ov{b(\varrho)},
\eF 
where the quantity $p - ({4 \over 3} \nu + \eta) \Div (\vc{u})$ is called the effective viscous pressure. 

In order to establish the same relation for variable viscosity coefficients, we write
\[
{\cal R}: [\psi \tn{S}] = \psi \left({4 \over 3} \nu + \eta \right) \Div (\vc{u}) + 
\Big\{ {\cal R} : [\psi \tn{S}] - \psi ({4 \over 3} \nu + \eta ) \Div (\vc{u}) \Big\}.
\]
Applying Lemma 4.2 in \cite{EF71}, the quantity
\[
{\cal R}: \Big[ \psi \Big( \nu (\tD, \bD) \left\langle \Grad \uD \right\rangle + \eta (\tD, \bD) \Div (\uD) \tn{I} \Big) \Big] - 
\psi ({4 \over 3} \nu (\tD, \bD) + \eta (\tD, \bD)) \Div (\uD),
\]
is bounded in the space $L^2(0,T; W^{\omega,p}(\Omega))$ for suitable $0 < \omega < 1$, $p > 1$ in terms of the 
bounds established in \eqref{u16}, \eqref{u17} and \eqref{u24} provided that the viscosity coefficients are 
globally Lipschitz functions of their arguments.

 If this is the case, one can use \eqref{e10} in order to deduce the 
desired relation  
\[
\Ov{ {\cal R}: [\psi \ \tn{S}] b(\varrho) } - \Ov{ {\cal R}: [\psi \ \tn{S}] } \ \Ov{b(\varrho)} = 
\Ov{ \psi ({4 \over 3} \nu + \eta )\Div (\vc{u}) \ b(\varrho)} - 
\Ov{ \psi ({4 \over 3} \nu + \eta) \Div(\vc{u}) } \ \Ov{ b(\varrho)} =
\]
\[
\psi ({4 \over 3} \nu + \eta )\Ov{ \Div(\vc{u}) \ b(\varrho)} - 
\psi ({4 \over 3} \nu + \eta) \Div(\vc{u})  \ \Ov{ b(\varrho)},
\]
where the last equality is a direct consequence of the pointwise convergence proved in \eqref{s8} and \eqref{e13}.

If $\nu$ and $\eta$ are only continuously differentiable as required by hypotheses of Theorem \ref{Tm1}, one can write
\bFormula{d4}
\left\{ 
\begin{array}{c}
\nu (\vartheta, \vc{B}) = Y (\vartheta, \vc{B}) \nu (\vartheta, \vc{B}) + (1 - Y(\vartheta, \vc{B})) 
\nu (\vartheta, \vc{B}), \\ \\
\eta (\vartheta, \vc{B}) = Y (\vartheta, \vc{B}) \eta (\vartheta, \vc{B}) + (1 - Y(\vartheta, \vc{B})) 
\eta (\vartheta, \vc{B}),
\end{array}
\right.
\eF
where $Y \in \D(\R^2)$ is a suitable function. 

Now we have
\[
\nu( \tD , \bD) \left\langle \uD \right\rangle + \eta (\tD, \bD) \Div(\uD) \tn{I} =
\]
\[
Y(\tD, \bD) \nu( \tD , \bD) \left\langle \uD \right\rangle +  Y(\tD, \bD) \eta (\tD, \bD) \Div(\uD) \tn{I} +
\]
\[
\Big\{ (1 - Y(\tD, \bD)) \nu( \tD , \bD) \left\langle \uD \right\rangle + (1 - Y(\tD, \bD)) \eta (\tD, \bD) \Div(\uD) \tn{I} \Big\},
\]
where the expression in the curl brackets can be made arbitrarily small in the norm of 
$L^p((0,T) \times \Omega)$, with a certain $p > 1$, by a suitable choice of $Y$ (see estimates 
\eqref{u5}, \eqref{u12}, \eqref{u17} and formula \eqref{u25}.  

Thus relation \eqref{d3} holds for any $\nu$ and $\eta$ satisfying hypotheses of Theorem \ref{Tm1} 
\vskip0,25cm
{\it 3. The oscillations defect measure}. In order to describe oscillations in the sequence $\sD{ \rD }$, we use the renormalized continuity equation 
\eqref{v3} together with its counterpart resulting from letting $\delta \to 0$ in \eqref{e9}.

 Although we have already shown that 
the limit quantities $\varrho$ and $\vc{u}$ satisfy the momentum equation \eqref{v1}, the validity of its renormalization \eqref{v3}
is not clear as the sequence $\sD{ \rD }$ is not known to be uniformly square integrable and the Di Perna-Lions machinery  
\cite{DL} does not work a priori. 

In order to solve this problem, a concept of {\it oscillations defect measure} was introduced in 
\cite{EF54}.

 To be more specific we set 
\bFormula{d5}
{\bf osc}_q[\rD \to \varrho] (Q) = \sup_{k \geq 1} \Big( \limsup_{\delta \to 0} 
\int_Q | T_k(\rD) - T_k(\varrho) |^q \ \dxdt \Big) ,
\eF
where $T_k$ are the cut-off functions 
\[
T_k(z) = {\rm sgn}(z) \min \{ |z| , k \}.
\] 

As shown in \cite{EF54}, the limit functions $\varrho$, $\vc{u}$ solve the renormalized equation \eqref{v3} provided
\begin{itemize}
\item
$\rD$, $\uD$ satisfy \eqref{e9}, 
\item
$\sD{ \uD }$ is bounded in $L^2(0,T; W^{1,2}(\Omega))$,
\item
\bFormula{d6}
{\bf osc}_q [ \rD \to \varrho] ((0,T) \times \Omega) < \infty \ \mbox{for a certain}\ q > 2.
\eF
\end{itemize}
Then, in order to establish \eqref{v3} it is clearly enough to show \eqref{d6}.
\vskip0.25cm
{\it 4. Weighted estimates of the oscillations defect measure}. In order to prove \eqref{d6} we use weighted estimates,
 where the corresponding weight function depends on the absolute temperature. 

Taking $b = T_k$ in \eqref{d3} we get
\bFormula{d7}
\Ov{p (\varrho, \vartheta) T_k(\varrho)} - \Ov{p(\varrho, \vartheta)}\ \Ov{T_k(\varrho)} 
\eF
\[
=\Big( {4 \over 3} \nu(\vartheta, \vc{B}) + \eta (\vartheta, \vc{B}) \Big) \Ov{ \Div(\vc{u}) \  T_k (\varrho) } - 
\Big( {4 \over 3} \nu(\vartheta, \vc{B}) + \eta(\vartheta, \vc{B}) \Big) \Div(\vc{u}) \ \Ov{ T_k (\varrho)}.
\] 
As observed in \eqref{p8}, there is a positive constant $p_c$ such that $p_F(\varrho, \vartheta) - p_c \varrho^{5 \over 3}$
is a non-decreasing function of $\varrho$ for any $\vartheta$.

 Accordingly we have 
\bFormula{d8}
\Ov{p (\varrho, \vartheta) T_k(\varrho)} - \Ov{p(\varrho, \vartheta)}\ \Ov{T_k(\varrho)} \geq
p_c \Big( \Ov{ \varrho^{5 \over 3} T_k (\varrho) } - 
\Ov{ \varrho^{5 \over 3} } \ \Ov{T_k (\varrho) } \Big).
\eF

Now let us choose a weight function $w \in C^1[0, \infty)$ 
\bFormula{d9}
w (\vartheta) > 0 \ \mbox{for}\ \vartheta \geq 0,\ 
w (\vartheta) = \vartheta^{-{1 + \alpha \over 2}} \ \mbox{for}\ \vartheta \geq 1,
\eF
where $\alpha $ is the exponent appearing in hypothesis \eqref{b31}.

Multiplying \eqref{d7} by $w(\vartheta)$ and using \eqref{d8} we obtain
\bFormula{d10} 
{ w(\vartheta)} \Big( 
 \Ov{ \varrho^{5 \over 3} T_k (\varrho) } - 
\Ov{ \varrho^{5 \over 3} } \ \Ov{T_k (\varrho) } \Big) 
\eF
\[
\leq { w(\vartheta) \over p_c }
\Big( {4 \over 3} \nu(\vartheta, \vc{B}) + \eta (\vartheta, \vc{B}) \Big) \Ov{ \Div(\vc{u}) \  T_k (\varrho) } - 
{ w(\vartheta) \over p_c} \Big( {4 \over 3} \nu(\vartheta, \vc{B}) + \eta(\vartheta, \vc{B}) \Big) \Div(\vc{u}) \ \Ov{ T_k (\varrho)},
\]
where the right-hand side is a weak limit of 
\[
{ w(\tD) \over p_c} \Big( {4 \over 3} \nu (\tD, \bD) + \eta (\tD, \bD) \Big) \Div(\uD) \Big) \Big( T_k (\rD) - T_k(\varrho) \Big)
+ 
{w(\vartheta) \over p_c} \Big( {4 \over 3} \nu(\vartheta, \vc{B}) + \eta(\vartheta, \vc{B}) \Big) \Div(\vc{u}) \Big( T_k(\varrho) 
- \Ov{ T_k (\varrho)} \Big).
\] 
Consequently, using the uniform weighted estimates \eqref{u12} together with hypothesis \eqref{b13} and the growth 
restriction  on $w$ specified in \eqref{d9}, we conclude that 
\bFormula{d11} 
\int_0^T \intO{ \Big[ 
{ w(\vartheta) \over p_c }
\Big( {4 \over 3} \nu(\vartheta, \vc{B}) + \eta (\vartheta, \vc{B}) \Big) \Ov{ \Div(\vc{u}) \  T_k (\varrho) } - 
{ w(\vartheta) \over p_c} \Big( {4 \over 3} \nu(\vartheta, \vc{B}) + \eta(\vartheta, \vc{B}) \Big) \Div(\vc{u}) \ \Ov{ T_k (\varrho)}
\Big]  } \ \dt
\eF
\[
\leq c \ \liminf_{\delta \to 0} \| T_k(\rD) - T_k (\varrho) \|_{L^2((0,T) \times \Omega)},\ \mbox{with}\ 
c \ \mbox{independent of}\ k .
\]
On the other hand, it was shown in Chapter 6 in \cite{EF70} that the left-hand side of \eqref{d10} could be bounded below as 
\bFormula{d12}
\int_0^T \intO{ 
{ w(\vartheta)} \Big( 
 \Ov{ \varrho^{5 \over 3} T_k (\varrho) } - 
\Ov{ \varrho^{5 \over 3} } \ \Ov{T_k (\varrho) } \Big) } \ \dt \geq
\limsup_{\delta \to 0} \int_0^T \intO{ w(\vartheta) | T_k(\rD) - T_k(\varrho) |^{8 \over 3} } \ \dt.
\eF
Then \eqref{d10} together with \eqref{d11}, \eqref{d12}, give the weighted estimate
\bFormula{d13}
\limsup_{\delta \to 0} \int_0^T \intO{ w(\vartheta) | T_k(\rD) - T_k(\varrho) |^{8 \over 3} } \ \dt \leq
c \ \liminf_{\delta \to 0} \| T_k(\rD) - T_k (\varrho) \|_{L^2((0,T) \times \Omega) }. 
\eF
Now taking $q > 2$ and $ {8 \over 3q} \omega  = {1 + \alpha \over 2}$, we can use H\"older's inequality to obtain 
\bFormula{d14}
\int_0^T \intO{ | T_k (\rD ) - T_k(\varrho) |^q } \ \dt
= \int_0^T \intO{ | T_k (\rD ) - T_k(\varrho) |^q (1 + \vartheta)^{-\omega} (1 + \vartheta)^{\omega} } \ \dt 
\eF
\[
\leq c \Big( \int_0^T \intO{ |T_k(\rD ) - T_k(\varrho) |^{8 \over 3} (1 + \vartheta)^{-{1 + \alpha \over 2}} } \ \dt
+ \int_0^T \intO{ (1 + \vartheta)^{ 3 q ( 1 + \alpha) \over 2 ( 8 - 3q) } } \ \dt \Big).
\]
 On the other hand, in accordance with \eqref{u5}, \eqref{u17} $\vartheta \in L^{\infty}(0,T; L^4(\Omega)) \cap L^3(0,T; L^9(\Omega))$, 
and a simple interpolation argument yields $\vartheta \in L^r((0,T) \times \Omega), \ \mbox{with}\ r = {46 \over 9}$. 
Consequently, one can find $q > 2$ such that 
\[
{3q (1 + \alpha) \over 2 ( 8 - 3q )} = r = {46 \over 9},
\]
provided $\alpha$ complies with hypothesis \eqref{m1}.

 Then finally relations \eqref{d13}, \eqref{d14} give rise to the desired 
estimate
\[
\limsup_{\delta \to 0} \int_0^T \intO{ w(\vartheta) | T_k(\rD) - T_k(\varrho) |^{q} } \ \dt \leq 
c \Big( 1 +  \liminf_{\delta \to 0} \| T_k(\rD) - T_k (\varrho) \|_{L^2((0,T) \times \Omega)} \Big),
\]
yielding \eqref{d6}.
\vskip0.25cm
{\it 5. Propagation of oscillations and strong convergence}. According to \eqref{d6} we know that the renormalized continuity equation \eqref{v3} is satisfied by the limit functions 
$\varrho$ and $\vc{u}$. In particular
\bFormula{d15}
\partial_t L_k(\varrho) + \Div ( L_k(\varrho) \vc{u} ) + T_k(\varrho) \ \Div(\vc{u}) = 0 
\ \mbox{in}\ \D'((0,T) \times \R^3),
\eF
provided $\varrho$, $\vc{u}$ have been extended to be zero outside $\Omega$,
where $\displaystyle L_k(\varrho) = \varrho \int_1^{\varrho} {T_k(z) \over z^2 } \ {\rm d}z$ .

On the other hand, one can let $\delta \to 0$ in \eqref{e9} and one obtains 
\bFormula{d16}
\partial_t \Ov{L_k(\varrho)} + \Div ( \Ov{L_k(\varrho)} \vc{u} ) + \Ov{ T_k(\varrho) \ \Div(\vc{u})} = 0 
\ \mbox{in}\ \D'((0,T) \times \R^3).
\eF
 Now, following step by step Chapter 3 in \cite{FEINOV} we take the difference of \eqref{d15} and \eqref{d16} and use 
\eqref{d7} and \eqref{d8} to deduce
\[
\intO{ \Ov{ L_k (\varrho) }(\tau) - L_k(\varrho)(\tau) } \leq \int_0^{\tau} 
\intO{ \left(\Div(\vc{u})\right) \Big( T_k(\varrho) - \Ov{T_k (\varrho) } \Big) } \ \dt 
\ \mbox{for any}\ \tau \in [0,T],
\]
where the right-hand sides vanishes for $k \to \infty$ due to \eqref{d6}.

Consequently 
\[
\Ov{ \varrho \log(\varrho) } = \varrho \log(\varrho)\ \ \ \mbox{in}\ (0,T) \times \Omega, 
\]
a relation equivalent to strong convergence of $\sD{ \rD }$, which means in fact that
\bFormula{d17}
\rD \to \varrho \ \mbox{in}\ L^1((0,T) \times \Omega).
\eF

\subsubsection{End of the proof}
\label{c}

 Note that we have already shown that $\varrho$, 
$\vc{u}$ satisfy the continuity equation \eqref{v1} as well as its renormalized version \eqref{v3}. 
Moreover, it is easy to see that the ``averaged'' momentum equation \eqref{s9} coincides in fact with \eqref{v4} in view of the strong 
convergence results established in \eqref{s8}, \eqref{e13} and \eqref{d17}. 
By the same token, one can pass to the limit in the energy equality \eqref{a29} in order to obtain \eqref{v12}. Note that 
the regularizing $\delta-$dependent terms on the left-hand side disappear by virtue of the estimates \eqref{u4}, 
\eqref{u5}, and \eqref{u28}.

Similarly, one can handle Maxwell's equation \eqref{a28}. Here, the only thing to observe is that the terms
\[
{1\over\sigma (\rD, \tD, \bH )} \Curl \bH \ \mbox{are bounded in}\ L^p((0,T) \times \Omega), \ 
\mbox{for a certain}\ p > 1,
\]
uniformly with respect to $\delta$ (such a bound can be obtained exactly as in \eqref{u25}.

To conclude, we have to deal with the entropy inequality \eqref{a27}. First of all, the extra terms on the left-hand side 
tend to zero because of \eqref{u4}, \eqref{u5}, \eqref{u22} and \eqref{e4}. Furthermore, it is standard to pass to 
the limit in the entropy production rate keeping the correct sense of the inequality as all terms are convex with respect 
to the spatial gradients of $\vc{u}$, $\vartheta$ and $\vc{B}$.

 Finally, the limit in the logarithmic terms 
can be carried over thanks to estimate \eqref{u22} and the following result (see Lemma 5.4 in \cite{DUFE1} holds
\bLemma{c1}
Let $\Omega \subset \R^N$ be a bounded Lipschitz domain. Assume that 
\[
\tD \to \vartheta \ \mbox{in}\ L^2((0,T) \times \Omega) \ \mbox{and}\ 
\log( \tD ) \to \Ov{log(\vartheta)} \ \mbox{weakly in}\ L^2((0,T) \times \Omega).
\]

Then $\vartheta$ is positive a.e. on $(0,T) \times \Omega$ and $\log(\vartheta) = \Ov{\log(\vartheta)}$.
\eL  

\section{Other models}
\label{other}
\subsection{Polytropic models}

The magnetic Navier-Stokes-Poisson case where magnetic permeability $\mu$ is a pure positive constant was studied by Ducomet and Feireisl in \cite{DUFE2}.

In the absence of magnetic field a model of  pressure law more adapted to cold plasmas of the type
\bFormula{b21bis}
p(\varrho, \vartheta) = p_e(\varrho) +  \vartheta p_{\vartheta}(\varrho),
\eF
with associated internal energy $e(\varrho, \vartheta)$ such that
\bFormula{b22bis}
e_{\rho}=\frac{1}{\rho^2}(p-\theta p_{\theta}),\ \ e_{\vartheta}=C_V(\vartheta),
\eF
introduced by Feireisl \cite{EF70} in the case of constant viscosities and 
 generalized to temperature dependent viscosities and conductivity in \cite{EF71}, has been 
extended to the MHD case (with Dirichlet boundary condition $\left. \vc{B}\right|_{\partial \Omega}=0$ and constant resistivity)
by Hu and Wang in \cite{HW1}:
\bTheorem{m1bis}
Let $\Omega \subset \R^3$ be a bounded domain with boundary of class $C^{2 + \delta}$, $\delta > 0$.

Suppose that the pressure $p=p(\varrho, \vartheta)$ is given by \eqref{b21bis}, that $p_e,p_{\vartheta}$ are $C^1$ functions of their arguments
on $[0,\infty)$ such that
\[
p_e(0)=p_{\vartheta}(0)=0,
\]
\[ 
p_e'(\rho)\geq a_1\rho^{\gamma-1},\ \ p_{\vartheta}(\rho)\geq 0\ \ \mbox{for all}\ \rho>0,
\]
\[
p_e(\rho)\leq a_2\rho^{\gamma},\ \ p_{\vartheta}(\rho)\leq a_3(1+\rho^{\gamma/3})\ \ \mbox{for all}\ \rho\geq 0.
\]
Furthermore, we suppose that the magnetic permeability $\mu$ and the electrical conductivity $\sigma$ are positive constants and that the transport coefficients 
\[
\kappa = \kappa( \vartheta),\ 
\nu = \nu(\vartheta),\ 
\frac 1 {\sigma\mu} := \lambda = \lambda(\vartheta):=\eta(\vartheta)-2/3,
\ c_V(\theta),
\]
are continuously differentiable functions of their argument satisfying 
\bFormula{B31bis}
c_1 (1 + \vartheta^{\alpha}) \leq \kappa( \vartheta) \leq c_2 (1 + \vartheta^{\alpha}),
\eF
\bFormula{B32bis}
0 \leq \lambda(\vartheta) \leq \overline{\lambda},
\eF 
\bFormula{B33bis}
0 <\underline{\nu}\leq \nu(\vartheta) \leq \overline{\nu},
\eF 
and
\bFormula{B34bis}
0 <\underline{c_V}\leq c_V(\vartheta) \leq \overline{c_V},
\eF 
for $ c_1,\ c_2,\ \overline{\lambda},\ \underline{\nu},\ \overline{\nu},\ \underline{c_V}, \overline{c_V}>0$ and $\alpha>2$.

Finally assume the homogeneous boundary conditions
\[
\left. \vc{H}\right|_{\partial\Omega}  = 0,\ \ 
\left.\vc{u}\right|_{\partial\Omega} = 0,\ \
\vc{q} \cdot \vc{n}|_{\partial \Omega} = 0,
\]
and let initial data $\varrho_0$, $(\varrho \vc{u})_0$, $\vartheta_0$, $\vc{B}_0$ be given so that
\bFormula{m2bis}
\varrho_0 \in L^{5 \over 3}(\Omega),\ 
(\varrho \vc{u})_0 \in L^1(\Omega; \R^3),\ \vartheta_0 \in L^{\infty} (\R^3),\ \vc{B}_0 \in L^2(\Omega; \R^3),
 \eF
\bFormula{m3bis}
\varrho_0 \geq 0,\ \vartheta_0  > 0, 
\eF
\bFormula{m4bis}
(\varrho s)_0 = \varrho_0 s(\varrho_0, \vartheta_0),\ 
{1 \over \varrho_0 } |(\varrho \vc{u})_0 |^2 ,\ (\varrho e)_0 = \varrho_0 e(\varrho_0, \vartheta_0) 
\in L^1(\Omega),
\eF
and
\bFormula{m5bis}
\Div \ \vc{B}_0 = 0 \ \mbox{in}\ \D'(\Omega),\ \vc{B}_0 \cdot \vc{n}|_{\partial \Omega} = 0.
\eF
Then problem \eqref{v1}-\eqref{v10} possesses at least one variational solution $\varrho$, $\vc{u}$, 
$\vartheta$, $\vc{B}$ on an arbitrary time interval $(0,T)$.
\eT

This model, with boundary conditions
$$\left. \vc{B}\cdot\vc{n}\right|_{\partial \Omega}=0,\ \ \left.  \Curl \vc{B}\times\vc{n}\right|_{\partial \Omega}=0,$$
 has been extended to the variable conductivity case
\[
0<\underline{\sigma}\leq \sigma(\varrho, \vartheta)\leq \overline{\sigma},
\]
for constants $\underline{\sigma},\ \overline{\sigma}>0$ and for $\alpha\geq 2$, but under additional constraints on viscosities, by Fan and Yu \cite{FY}.

\subsection{Barotropic models}

When the temperature fluctuations can be neglected, one obtains the barotropic-MHD system for the density $\varrho$, the velocity $\vc{u}$
and the magnetic field $\vc{H}$ in a bounded region $\Omega$
\bFormula{baro1}
\partial_t \varrho + \Div (\varrho \vc{u}) = 0,
\eF
\bFormula{baro2}
\partial_t (\varrho \vc{u}) + \Div (\varrho \vc{u} \otimes \vc{u}) = \Div \tn{T} + \vc{j} \times \vc{B},
\eF 
\bFormula{baro3}
\partial_t \vc{H} + \Curl( \vc{H} \times \vc{u} ) + \Curl( \lambda \Curl  \vc{H} ) = 0,\ \ \ \ \Div \vc{H} = 0,
\eF
\bFormula{baro4}
 \left. \vc{u}\right|_{\partial \Omega} = 0 ,\ \left. \vc{H} \right|_{\partial \Omega} = 0.
\eF
where the pressure is $p(\varrho)=A\varrho^{\gamma}$ and the transport coefficients $\nu,\eta,\lambda$ are positive constants.

Hu and Wang proved the following existence result \cite{HW1}
\bTheorem{HW}
Let $\Omega \subset \R^3$ be a bounded domain with boundary of class $C^{2 + \delta}$, $\delta > 0$ and $\gamma>3/2$.
Assume that the initial data $\varrho_0$, $(\varrho \vc{u})_0$ and $\vc{B}_0$ be given so that
\bFormula{mm1}
\varrho_0 \in L^{\gamma}(\Omega),\ 
(\varrho \vc{u})_0 \in L^1(\Omega; \R^3),\ \vc{H}_0 \in L^2(\Omega; \R^3),
 \eF
\bFormula{mm2}
\varrho_0 \geq 0,
\eF
and
\bFormula{mm3}
\Div  \vc{H}_0 = 0 \ \mbox{in}\ \D'(\Omega),\ \vc{H}_0 |_{\partial \Omega} = 0.
\eF
Assume finally that $f,g\in L^{\infty}(0,T\times\Omega)$.

Then problem \eqref{baro1}-\eqref{baro4} possesses at least one variational solution $\varrho$, $\vc{u}$, $\vc{H}$ on an arbitrary time interval $(0,T)$.
\eT
Let us mention that the same result holds (R. Sart \cite{S}) with different boundary conditions on $\vc{H}$.
\vskip0.25cm
Let us quote for completeness some other models recently introduced.

 Y. Amirat and K. Hamdache \cite{AHa} studied a barotropic ferrofluid model with a 
magnetization dependent pressure $p=p(\varrho,\vc{M})$ including in the fluid description an extra equation for the angular momentum.

Finally X. Hu and D. Wang \cite{HW3} considered density dependent viscosities, extending to MHD the result of D. Bresch and B. Desjardins \cite{BD} relying 
on analysis of the so-called BD-entropy. In this last case, presence of the radiative Stefan-Boltzmann contribution $\vartheta^4$ is no-more requested but
additional constraints on the pressure and transport coefficients (including magnetic permeability) are needed.

\section{Singular limits}
\label{singular}

As usual in hydrodynamical problems, when putting the system in non dimensional form, various non dimensional numbers \cite{ZE} (Mach, Reynolds, Peclet etc...)
appear, which may be small (or large) in certain regimes of motion. Studying the corresponding limits may simplify drastically the description. 
Typically low Mach number limits of compressible flows are expected to be close to their incompressible limit.
In MHD flows, asymptotic regimes appear naturally when the so called Mach number, Alfven number, Froude number or P\'eclet number are small.

The simplest case has been studied by  P. Kuku\v cka \cite{K}, for which one expects convergence toward an incompressible limit.
\vskip0.25cm
More precisely for a suitable scaling (see \cite{K}) one considers the following  non dimensional form of the MHD system \eqref{b7}-\eqref{b8}-\eqref{b9}-\eqref{b16}
where magnetic permeability $\mu$ is supposed to be a positive constant 
(and we note $\lambda(\varrho,\vartheta,\vc{H}):=\frac{1}{\mu\sigma(\varrho,\vartheta,\vc{H})}$ the magnetic
diffusivity), to simplify the exposition.
\bFormula{bma7}
\partial_t \vc{B} + \Curl( \vc{B} \times \vc{u} ) + \Curl\left( \lambda \Curl  \vc{B}\right) = 0,\
\eF
\bFormula{bma8}
\partial_t \varrho + \Div (\varrho \vc{u}) = 0,
\eF
\bFormula{bma9}
\partial_t (\varrho \vc{u}) + \Div (\varrho \vc{u} \otimes \vc{u})\frac{1}{Ma^2} = \Div \tn{S} +\frac{1}{Al^2}\vc{j} \times \vc{B}+\frac{1}{Fr^2}\varrho \Grad \Phi,
\eF 
\bFormula{bma10}
\frac{d}{dt} \Big( {Ma^2 \over 2} \varrho |\vc{u}|^2 + \varrho e  +\frac{Ma^2}{Al^2}\frac{1}{2\mu}|\vc{B}|^2\Big)\ dx=0,
\eF
\bFormula{bma11}
\partial_t (\varrho s) + \Div (\varrho s \vc{u} ) + \frac{1}{Pe}\Div \Big({\vc{q} \over \vartheta } \Big) = \sigma,
\eF
with
\bFormula{bma12}
\sigma \geq {1 \over \vartheta}
\Big( Ma^2\tn{S} : \Grad \vc{u} + \frac{Ma^2}{Al^2} |\Curl \vc{H}|^2 - \frac{1}{Pe}{\vc{q} \cdot \Grad \vartheta \over \vartheta} \Big),
\eF
where $\Phi$ is an external potential (gravitation for example).

In order to simplify the exposition, following \cite{K}, we suppose that Mach number $Ma$ and Alfven number $Al$ are both small 
$Ma=Al=\varepsilon$ with $0<\varepsilon\ll 1$ while $Pe=1$. 
We also discard the external potential force ($\Phi=0$).

We assume boundary conditions
\bFormula{bcbma}
\left. \vc{B}\cdot \vc{n}\right|_{\partial\Omega}  = 0,
\ \ \left. \vc{E}\times \vc{n}\right|_{\partial\Omega}  = \vc{0},
\ \ \left. \vc{u}\cdot \vc{n}\right|_{\partial\Omega}  = 0,
\ \ \left. \tn{S}\vc{n}\times \vc{n}\right|_{\partial\Omega}  = \vc{0},
\ \ \left. \Grad \vartheta\cdot \vc{n}\right|_{\partial\Omega}  = \vc{0},
\eF
and initial conditions
\bFormula{icbma}
(\vc{B},\varrho,\vc{u},\vartheta)(0,x)
 =(\vc{B}_0,\varrho_0,\vc{u}_0,\vartheta_0)(x)\ \ \mbox{for any}\ x\in\Omega,
\eF
``ill-prepared" in the sense that
\bFormula{iceps}
\ \vc{B}_0(x)= \vc{B}_{0,\varepsilon}^{(1)}(x),
\ \varrho_0(x)=\overline{\varrho}+\varepsilon_{0,\varepsilon}^{(1)}(x),
\ \vc{u}_0(x)= \vc{u}_{0,\varepsilon}^{(1)}(x),
\ \vartheta_0(x)=\overline{\vartheta}+\vartheta_{0,\varepsilon}^{(1)}(x),
\eF
where $\overline{\varrho}$ and $\overline{\varrho}$ are two positive constants.

Formally expanding \eqref{bma7}-\eqref{bma11} with respect to $\varepsilon$
 we actually find at lowest order the incompressible MHD system with temperature
\bFormula{ba7}
\partial_t \overline{\vc{B}} + \Curl( \overline{\vc{B}} \times \vc{U} ) +
 \Curl\left( \lambda(\overline{\varrho},\overline{\vartheta},\overline{\vc{B}}) \Curl  \overline{\vc{B}}\right) = 0,
\eF
\bFormula{ba8}
 \Div (\vc{U}) = 0,
\eF
\bFormula{ba9}
\overline{\varrho}\left(\partial_t \vc{U} +  \vc{U} \cdot\Grad \vc{U}\right)+\Grad \Pi
 = \Div \left(\tn{S}(\vc{U},\overline{\vartheta})\right) +\Curl \overline{\vc{B}} \times \overline{\vc{B}},
\eF 
\bFormula{ba11}
\overline{\varrho}c_P(\overline{\varrho},\overline{\vartheta})\left(\partial_t \Theta +  \vc{U} \cdot\Grad \Theta\right)
 = \Div \left(\left(\kappa_F(\overline{\varrho},\overline{\vartheta},\overline{\vc{B}})+\kappa_r \overline{\vartheta}^2\right)\Grad \Theta\right),
\eF
where $c_P(\varrho,\vartheta)=e_{\vartheta}+\frac{\vartheta}{\varrho^2}\frac{p_{\vartheta}^2}{p_{\varrho}}$, with boundary conditions
\bFormula{bcbmat}
\left. \vc{B}\cdot \vc{n}\right|_{\partial\Omega}  = 0,
\ \ \left. \vc{U}\cdot \vc{n}\right|_{\partial\Omega}  = 0,
  \ \ \left. \Grad \Theta\cdot \vc{n}\right|_{\partial\Omega}  = 0,
\eF
and initial conditions
\bFormula{icbmat}
\vc{B}(0,x)=\vc{B}_0(x),\ \vc{U}(0,x)=\vc{u}_0(x),\ \Theta(0,x)=\vartheta_0(x)\ \ \mbox{for any}\ x\in\Omega.
\eF
then one has the convergence result \cite{K}
\bTheorem{kU}
Let $\Omega \subset \R^3$ be a bounded domain with boundary of class $C^{2 + \delta}$, $\delta > 0$ and 
assume that properties of the state functions and transport coefficients given in Theorem \ref{Tm1} are valid.

Let $\{\varrho_{\varepsilon},\vc{u}_{\varepsilon},\vartheta_{\varepsilon},\vc{B}_{\varepsilon}\}_{\varepsilon>0}$ be
a family of weak solutions to the scaled MHD system \eqref{bma7}-\eqref{bma11} on $(0,T)\times\Omega$ supplemented with boundary conditions
\eqref{bcbma} and initial conditions \eqref{iceps}
such that
\[
\varrho_{0,\varepsilon}\to \varrho_0^{(1)}\ \ \mbox{weakly* in}\ L^{\infty}(\Omega),
\]
\[
 \vc{u}_{0,\varepsilon}\to \vc{u}_0^{(1)}\ \ \mbox{weakly* in}\ L^{\infty}(\Omega;\R^3),
\]
\[
\vartheta_{0,\varepsilon}\to \vartheta_0^{(1)}\ \ \mbox{weakly* in}\ L^{\infty}(\Omega),
\]
\[
 \vc{B}_{0,\varepsilon}\to \vc{B}_0^{(1)}\ \ \mbox{weakly* in}\ L^{\infty}(\Omega;\R^3),
\]
as $\varepsilon\to 0$. Then
\[
\mbox{ess sup}_{t\in(0,T)}\| \varrho_{\varepsilon}(t)-\overline{\varrho}\|_{L^{\frac{5}{3}}(\Omega)}\leq C\varepsilon,
\]
\[
 \vc{u}_{\varepsilon}\to \vc{U}\ \ \mbox{weakly in}\ L^2(0,T;W^{1,2}(\Omega;\R^3),
\]
\[
 \frac{\vc{B}_{\varepsilon}}{\varepsilon}\to \vc{B}\ \ \mbox{weakly in}\ L^2(0,T;W^{1,2}(\Omega;\R^3),
\]
\[
 \frac{\vartheta_{\varepsilon}-\overline{\vartheta}}{\varepsilon}\to \Theta \ \mbox{weakly in}\ L^2(0,T;W^{1,2}(\Omega;\R^3),
\]
where $(\vc{U},\Theta,\vc{B})$ is a weak solution of the incompressible MHD system with temperature \eqref{ba7}-\eqref{ba11}
 with boundary conditions \eqref{icbmat}, initial temperature
\[
\Theta_0=\frac{\overline{\vartheta}}{c_P(\overline{\varrho},\overline{\vartheta})}
\left(s_{\varrho}(\overline{\varrho},\overline{\vartheta})\varrho_0^{(1)}+s_{\vartheta}(\overline{\varrho},\overline{\vartheta})\vartheta_0^{(1)}\right),
\]
and initial velocity given by
\[
\overline{\varrho}\vc{U}(0)=\lim_{\varepsilon\to 0}H[  \varrho_{\varepsilon} \vc{u}_{\varepsilon}](0),
\]
where $H$ is Helmholtz's projector on the space of solenoidal fields.
\eT
The proof follows the stategy of \cite{FEINOV}:
\begin{enumerate} 
\item One first proves uniform estimates for the sequence $\{\varrho_{\varepsilon},\vc{u}_{\varepsilon},\vartheta_{\varepsilon},\vc{B}_{\varepsilon}\}_{\varepsilon>0}$
and the entropy rate $\sigma_{\varepsilon}$.
\item One pass to the limit $\varepsilon\to 0$ in the system, the most delicate term being the convective term $\varrho_{\varepsilon}\vc{u}_{\varepsilon}\otimes\vc{u}_{\varepsilon}$,
for which one uses the Helmholtz decomposition together with the associated acoustic equation.
\end{enumerate}

More complex situations have been investigated:
 A. Novotn\`y, M. R$\mathring{u}$\v zi\v cka ang G. Th\"ater \cite{NRT} have considered the case
$Ma= \epsilon$, $Al=\epsilon$, $Fr=\epsilon$, $Pe= \epsilon^2$, 
Y.S. Kwon and K. Trivisa \cite{KT} have studied the case
$Ma= \epsilon$, $Al=\sqrt{\epsilon}$, $Fr= \sqrt{\epsilon}$, $Pe=1$,
E. Feireisl,  A. Novotn\`y and Y. Sun \cite{FNS}  have considered 
$Ma= \epsilon$, $Al=Fr=Pe=1$ $Re=1/\epsilon^n$, for a $n>0$ (large Reynolds number regime) and
S. Jiang, Q. Ju, F. Li and Z. Xin \cite{JJLX} have studied more complex $\epsilon^k$ regimes. In all of these situations the target system
is the incompressible MHD system.

\section{Strong solutions}
\label{strong}

Previous Sections dealt with weak solutions, belonging to Sobolev spaces, for system
\eqref{b7}-\eqref{b8}-\eqref{b9}-\eqref{b14}, together with \eqref{O} and \eqref{Abis} for the definition of the
electric field and the current density. It is not known in general if these solutions are unique. 

Another approach is to look directly for strong solutions, which are assumed to be regular. In order to do so, the general
method is to prove local existence of smooth solutions. This is done by applying general local existence results on
hyperbolic and hyperbolic-parabolic systems (Subsection~\ref{strong_local} below). Then, it is possible, assuming
that the initial data are sufficiently close to an equilibrium state, to prove that the parabolic
nature of the equation dominates. This implies a priori bounds which allow to prove that the solution remains in
the neighbourhood of the equilibrium as long as it exists. This in turn allows to prove global existence. This
method, due to Kawashima and his co-authors, is described in Subsection~\ref{strong_global}. Subsection~\ref{critical_spaces} is devoted
to existence and uniqueness results in critical (Besov) spaces.

\subsection{Local existence results}
\label{strong_local}

In \cite{VH}, it is proved that a symmetric hyperbolic-parabolic system has a unique local solution. Therefore, its
application to the present system only amounts to proving that the system is symmetrizable. The use of entropy
variables is necessary to give a symmetric system for the fluid part (equations \eqref{b8}-\eqref{b9}-\eqref{b14}). See for instance \cite{D} for
the details. The part with the magnetic field is dealt with using the method of \cite{SE}, so the results of
\cite{VH} apply. The result may be summarized as follows:
\begin{Theorem}
   \label{th:strong_local} 
Let $\Omega = \R^3$, and assume that the thermodynamic functions
\[
p=p(\varrho, \vartheta), \ e = e(\varrho, \vartheta),\  
s = s(\varrho, \vartheta)
\]
are interrelated through \eqref{b17}, where $p$, $e$ are continuously differentiable functions of positive arguments $\varrho$, $\vartheta$ 
satisfying \eqref{b23}-\eqref{b24}-\eqref{b25}-\eqref{b26}. 

Furthermore, we suppose that each of the transport coefficients 
\[
\nu = \nu(\vartheta, |\vc{H}|),\ 
\eta = \eta(\vartheta, |\vc{H}|),\ \kappa_F = \kappa_F(\varrho, \vartheta, |\vc{H}|), \ \sigma = 
\sigma(\varrho, \vartheta, |\vc{H}|)\ \mbox{and} \ \mu=\mu( |\vc{H}|),
\]
are positive smooth functions of their arguments. Consider an equilibrium state $\overline
\varrho, \overline{\vc{u}} \equiv 0, \overline\vartheta, \overline{\vc{B}}\equiv 0$ of the system
\eqref{b7}-\eqref{b8}-\eqref{b9}-\eqref{b14}, and let the initial data data $\varrho_0$, $\vc{u}_0$, $\vartheta_0$, $\vc{B}_0$ be given so that
\begin{equation}
  \label{eq:DI_locale}
  \left(\varrho_0 - \overline\varrho, \vc{u}_0,\vartheta_0 - \overline{\vartheta}, \vc{B}_0 
  \right) \in H^{3}(\R^3).
\end{equation}
Then there exists $T>0$ such that the system \eqref{b7}-\eqref{b8}-\eqref{b9}-\eqref{b14}, with the initial
condition $\varrho_0$, $(\varrho \vc{u})_0$, $\vartheta_0$, $\vc{B}_0$, has a unique strong solution on
$[0,T]$. This solution satisfies $\varrho,\vartheta>0$ in $\R^3\times [0,T]$, and $T$ depends only on $\left\|
\left(\varrho_0-\overline\varrho,\vc{u}_0,\vartheta_0-\overline\vartheta,\vc{B}_0 \right)\right\|_{H^3}$, $\displaystyle\inf_{\R^3} \vartheta_0$ and $\displaystyle\inf_{\R^3} \varrho_0$.
\end{Theorem}

Some remarks are in order. 

First, the above result is stated in \cite{K2}, for electro-magneto-fluids (which, in
addition to the above systems, includes equations for the charge $\varrho_e$ and the electric field $\vc{E}$) with
two-dimensional symmetries. However, the aim of \cite{K2} is to prove global existence, and this is why their study
is restricted to two-dimensional flows. The local existence result is still valid for the 3D case. As a matter of
fact, it is an application of the results of \cite{VH}, which are valid in any dimension. 

Second, the same result can be stated with $H^s$ instead of $H^3$, for any $s\geq 3$ (see \cite[Section 4]{VH}).

Third, Theorem~\ref{th:strong_local} is proved in \cite[Theorem 2.1]{PG} in the special case of constant coefficients
$\nu,\eta,\kappa_F,\sigma$ and $\mu$. The proof is not given, and is claimed to be an easy adaptation of
that of \cite{MN}. Looking at the proof of \cite{MN}, it is restricted to the case $\vc{B} = 0$, that is, the
compressible Navier-Stokes system with thermal conductivity. Their proof readily applies to the present case, and
allows in fact for variable coefficients satisfying the above assumptions. A complete proof is also given in \cite{FY2},
in the case of a bounded domain (the method is easily adapted to $\Omega=\R^3$), with possible occurrence of vacuum,
that is, $\varrho\geq 0$ instead of $\varrho>0$.

Note also that, in Theorem~\ref{th:strong_local}, we have assumed \eqref{b23}, that is, the fluid is a perfect
gas. Now, the result is still valid for a general equation of state satisfying \eqref{b24}-\eqref{b25}-\eqref{b26}
only (see \cite{K2}).

Finally, the assumptions
$\nu,\eta,\kappa_F, \sigma, \mu >0$ may be weakened in the following way: it is in fact possible to treat cases in
which some of these coefficients are identically $0$. For instance, using the same method as in \cite{K2}, it is
possible to deal with the case $\nu,\eta, \sigma, \mu >0, \kappa_F \equiv 0$, the case $\nu,\kappa_F, \sigma, \mu
>0, \eta \equiv 0$, and the case $\nu, \sigma, \mu >0, \kappa_F, \eta \equiv 0$. 

\subsection{Global existence results for small data}
\label{strong_global}

Theorem~\ref{th:strong_local} being proved, one can then turn to the proof of global existence of a smooth solution
for the system under consideration. For this purpose, the key ingredient is an a priori estimate which we now
detail. For this purpose, we define, for any $0\leq t_1 \leq t_2$, the quantity
\begin{multline}
  \label{eq:N_s}
  N_s(t_1,t_2) = \sup_{t_1\leq t\leq t_2}
  \left\|\left(\varrho-\overline\varrho,\vc{u},\vartheta-\overline\vartheta, \vc{B} \right)(t) \right\|_{H^s(\R^3)}
  \\+ \left[\int_{t_1}^{t_2}\left(\left\|\left(\nabla\varrho,\nabla\vc{B}\right) \right\|_{H^{s-1}(\R^3)}^2(t) + \left\|
    \left(\nabla\vc{u},\nabla\vartheta\right)\right\|_{H^s(\R^3)}^2(t) \right)dt \right]^{1/2}.
\end{multline}
Then, differentiating the system if necessary, and multiplying it by some well chosen test functions, it is
possible to prove the following kind of estimate: there exists $\delta>0$ such that, if 
\begin{displaymath}
  N_s(0,t)^2 \leq C \left(N_s(0,0)^2 + N(0,t)^3\right),
\end{displaymath}
for some constant $C$ which depends on $\delta$, but not on $t>0$. This in turn allows to prove that there exists
$\delta>0$ and $C(\delta)>0$, both independent of $t$, such that
\begin{equation}
  \label{eq:estimation_kawashima}
N_3(0,0)\leq \delta \quad \Rightarrow \quad  N_3(0,t) \leq C(\delta) N_3(0,0).
\end{equation}
Using \eqref{eq:estimation_kawashima}, one proves a bound on $N_3(0,t)$ for any time $t$ of existence of the local in time solution. Applying
the local existence result successively any each small time interval, one finally proves global existence:
\begin{Theorem}
  \label{th:strong_global}
Let $\Omega = \R^3$, and assume that the thermodynamic functions
\[
p=p(\varrho, \vartheta), \ e = e(\varrho, \vartheta),\  
s = s(\varrho, \vartheta)
\]
are interrelated through \eqref{b17}, where $p$, $e$ are continuously differentiable functions of positive arguments $\varrho$, $\vartheta$ 
satisfying \eqref{b24}-\eqref{b26}. 

Furthermore, we suppose that each of the transport coefficients $\nu,\eta,\kappa_F,\sigma,\mu$ are positive constants.

Consider an equilibrium state $\overline
\varrho, \overline{\vc{u}} \equiv 0, \overline\vartheta, \overline{\vc{B}}\equiv 0$ of the system
\eqref{b7}-\eqref{b8}-\eqref{b9}-\eqref{b14}, and let the initial data $\varrho_0$, $\vc{u}_0$, $\vartheta_0$, $\vc{B}_0$ be given so that
\begin{equation}
  \label{eq:DI_locale_2}
  \left(\varrho_0 - \overline\varrho, \vc{u}_0,\vartheta_0 - \overline{\vartheta}, \vc{B}_0 
  \right) \in H^{3}(\R^3).
\end{equation}
Then there exists $\delta>0$ such that, if $N_3(0,0)<\delta$ ($N_3$ is defined by \eqref{eq:N_s}, then problem  \eqref{b7}-\eqref{b8}-\eqref{b9}-\eqref{b14}, with the initial
condition $\varrho_0$, $\vc{u}_0$, $\vartheta_0$, $\vc{B}_0$, has a unique global strong solution
$(\varrho,\vc{u},\vartheta, \vc{B})$ such that
\begin{multline}\label{eq:bornes_existence_globale}
  \left(\varrho-\overline{\varrho}, \vc{B}\right) in C^0\left(\R^+, H^3\right)\cap C^1\left(\R^+, H^2\right),\quad
  \left(\vc{u},\vartheta-\overline\vartheta \right) \in C^0\left(\R^+, H^3\right)\cap C^1\left(\R^+, H^1\right),
  \\ \nabla \left(\vc{u},\vartheta-\overline\vartheta \right)\in L^2\left(\R^+,H^3\right).
\end{multline}
This solution satisfies the estimate
\begin{displaymath}
\forall t\in \R^+, \quad  N_3(0,t) \leq C N_3(0,0),
\end{displaymath}
where $C$ depends only on $\delta>0$, and $N_3$ is defined by \eqref{eq:N_s}.
\end{Theorem}

This result is proved in \cite[Theorem 2.5]{PG}. Considering the strategy of proof, it seems that the fact that the
coefficients are constant is not important, although we did not check all the details. 

In \cite{K2}, the same result is proved for a more general system and with non-constant coefficients, but only with
two-dimensional symmetries. Note that, due to the structure of their system, which includes an equation for the
electric field and the density of charge, the equilibrium state does not necessarily satisfies $\overline{\vc{B}}
\equiv 0$. 

In essence, the above result is proved by studying the linearized equation near the equilibrium state
$\left(\overline\varrho, \overline{\vc{u}}\equiv 0, \overline\vartheta,\overline{\vc{B}}\equiv 0\right)$. This is a
linear hyperbolic-parabolic system, for which energy-type estimates allow to prove the above result. Hence, a
natural consequence is that the solution converges, as time tends to infinity, to the equilibrium state
$\left(\overline\varrho, \overline{\vc{u}}\equiv 0, \overline\vartheta,\overline{\vc{B}}\equiv 0\right)$. Such a
convergence has been studied in \cite{PG}, with the above assumptions. The convergence following convergence is
proved in this paper:
\begin{equation}
  \label{eq:decay1}
  \forall k=1,2,3, \quad \left\|\nabla^k \left(\varrho-\overline\varrho, \vc{u},\vartheta-\overline\vartheta, \vc{B}\right) \right\|_{L^2(\R^3)} \leq \frac C {(1+t)^{\frac32\left(\frac 1 q - \frac12 \right)+\frac k 2}}
\end{equation}
\begin{equation}
  \label{eq:decay2}
\forall p>q, \quad \left\|\left(\varrho-\overline\varrho, \vc{u},\vartheta-\overline\vartheta, \vc{B}\right)
\right\|_{L^p(\R^3)} \leq \frac C {(1+t)^{\frac32\left(\frac 1 q - \frac 1 p \right)}} 
\end{equation}
Here, $q\in \left[1,\frac65\right)$ is such that $\left(\varrho_0-\overline\varrho, \vc{u}_0,
  \vartheta_0-\overline\vartheta, \vc{B}_0\right) \in L^q(\R^3).$

\medskip

As in the case of variational solutions, barotropic fluids have been extensively studied. In such a
case, the system becomes \eqref{baro1}-\eqref{baro2}-\eqref{baro3}-\eqref{baro4}, with a pressure law 
$p(\varrho) = \varrho^\gamma,$ $\gamma>1$. All the coefficients are supposed to be constant. In \cite{LSW,WW}, local existence is proved in the case
$\lambda=0$, that is, when the parabolic term in \eqref{baro3} vanishes. The main idea of the method is to fix the
velocity field $u$, and then solve the problem with respect to all other variables. Then the Schauder fixed-point
theorem is applied. Note that in \cite{WW}, $\varrho\mapsto p(\varrho)$ is only assumed to be strictly convex and
increasing. Global existence is proved in \cite{CT} under the assumption of small initial data, using a method
similar to the one explained above. \cite{TTW} proposes a generalization to the case of a system with Coulomb
interaction, and \cite{YGD} to the case of a bounded domain. The result in \cite{LXZ} is similar, but allows for
vacuum in the system of equations. Finally, \cite{YGD} gives a similar kind of result, but with an $H^2$ norm to
estimate the initial data. Note finally that the asymptotic behavior as $t\to+\infty$ is also proved in
\cite{CT,TTW,YGD} with the same kind of results as above.

\subsection{Critical regularity}
\label{critical_spaces}

After the works of R. Danchin \cite{D1} \cite{D2} \cite{D3} on compressible Navier-Stokes system 
in a critical functional framework, some recent works have been 
devoted to the extension of his results to the compressible MHD system. 

For the reader's convenience, let us first recall a few basic definitions concerning 
Besov spaces (a much more detailed presentation can be found in \cite{BCD}). 
\medbreak
Let us first introduce the homogeneous Littlewood-Paley decomposition.
 Let  $\chi:\R^n\rightarrow[0,1]$ be  a smooth nonincreasing radial function
supported in $B(0,1)$ and such that $\chi\equiv1$ on $B(0,1/2),$ and let 
$$\varphi(\xi):=\chi(\xi/2)-\chi(\xi).$$ 
Then  we define the homogeneous Littlewood-Paley
decomposition  $(\dot\Delta_k)_{k\in\Z}$ over $\R^{n}$
by 
$$\dot\Delta_ku:=\varphi(2^{-k}D)u=
{\mathcal F}^{-1}\bigl(\varphi(2^{-k}\cdot){\mathcal F}u\bigr),
$$
where  $\cF$ stands for the Fourier transform on $\R^{n}.$

For any $s\in\R$ and $(p,r)\in[1,+\infty]^2,$ 
the homogeneous Besov space $\dot B^s_{p,r}(\R^n)$ stands 
for the set of tempered distributions $f$ such that 
$$
\|f\|_{\dot B^s_{p,r}(\R^{n})}:= \bigl\|2^{sk}\|\dot\Delta_kf\|_{L^p(\R^{n})}\bigr\|_{\ell^r(\Z)}<\infty,
$$
and we also note $B^s:=B^s_{2,2}$.

One has the formal decomposition
$$
u=\sum_{k\in\Z}\ddk u.
$$
Owing to the fact that the solutions that we consider here do not have the same regularity in low and high frequencies, 
it is useful to use hybrid Besov spaces by splitting any distribution  $u$ into its low frequency part and 
its high frequency part \cite{BCD} as follows:
\begin{equation}\label{eq:split}
u^\ell=\sum_{k\leq 0}\ddk u\quad\hbox{and}\quad u^h=\sum_{k>0} \ddk u.
\end{equation}
Then we define for $s,t\in\R$
\[
\|f\|_{B^{s,t}}:=\sum_{k\leq0}2^{sk}\|\ddk f\|_{L^2}+\sum_{k>0}2^{tk}\|\ddk f\|_{L^2},
\]
and the hybrid Besov space $B^{s,t}$ is defined for $m=-[d/2+1-s]$ by
\[
B^{s,t}:=\{f\in{\mathcal S}'(\R^d)\ :\ \|f\|_{B^{s,t}}<\infty\},
\]
if $m<0$, or
\[
B^{s,t}:=\{f\in{\mathcal S}'(\R^d)/{\mathcal P}_m\ :\ \|f\|_{B^{s,t}}<\infty\},
\]
if $m\geq 0$, where ${\mathcal P}_m$ is the set of polynomials of degree less or equal to $m$ (of course in any case $B^s:= B^{s,s}$).
\vskip0.25cm
Let us first rewrite the barotropic system
\bFormula{baroc1}
\partial_t \varrho + \Div (\varrho \vc{u}) = 0,
\eF
\bFormula{baroc2}
\partial_t (\varrho \vc{u}) + \Div (\varrho \vc{u} \otimes \vc{u}) -\Grad p(\varrho)= \mu\Delta\vc{u}+(\mu+\lambda)\Grad \Div(\vc{u}) +  \vc{H}\cdot\Grad \vc{H}-\frac{1}{2}\Grad(|\vc{H}|^2),
\eF 
\bFormula{baroc3}
\partial_t \vc{H} + \Div(\vc{u})\ \vc{H}+  \vc{u}\cdot \Grad \vc{H} ) =\nu \Delta \vc{H},\ \ \ \ \Div \vc{H} = 0,
\eF
\bFormula{baroc4}
 (\varrho, \vc{u},\vc{H})_{t=0} = (\varrho_0, \vc{u}_0,\vc{H}_0),\ \ \ \ \Div \vc{H}_0 = 0.
\eF
where the pressure is $p(\varrho)=A\varrho^{\gamma}$ and the transport coefficients $\nu,\eta,\lambda$ are positive constants.

Hao proved the following existence result \cite{Hao}
\bTheorem{HW1}
Let $d\geq 3$, $\overline{\rho}>0$, $2\mu+d\lambda>0$, $\nu>0$ and $\vc{h}$  an arbitrary constant vector.

Assume that the initial data satisfy $\varrho_0-\overline{\rho}\in B^{d/2-1,d/2}$, $ \vc{u})_0\in B^{d/2-1}$ and $\vc{H}_0\in B^{d/2-1}$.

Then there exists a small number $\varepsilon>0$ and a constant $M$ such that
\[
\|\varrho_0-\overline{\rho}\|_{ B^{d/2-1,d/2}}+\| \vc{u}_0\|_{ B^{d/2-1}}+\|\vc{H}_0\|_{ B^{d/2-1}}\leq \varepsilon.
\]
then system \eqref{baroc1}-\eqref{baroc4} has a unique global solution on an arbitrary time interval $(0,T)$ such that
$(\varrho-\overline{\rho},\vc{u},\vc{H}-\vc{h})\in E$, with
\[
E:=C\left(\R_+; B^{d/2-1,d/2}\times  (B^{d/2})^{d+d}\right)\cap L^1\left(\R_+; B^{d/2+1,d/2+2}\times  (B^{d/2+1})^{d+d}\right),
\]
and satisfies
\[
\|(\varrho-\overline{\rho},\vc{u},\vc{H}-\vc{h})\|_E
\leq M
\big(\|\varrho_0-\overline{\rho}\|_{ B^{d/2-1,d/2}}+\| \vc{u}_0\|_{ B^{d/2-1}}+\|\vc{H}_0\|_{ B^{d/2-1}}\big).
\]
\eT
In the case of a polytropic fluid with constant viscosities the system rewrites
\bFormula{barop1}
\partial_t \varrho + \Div (\varrho \vc{u}) = 0,
\eF
\bFormula{barop2}
\partial_t (\varrho \vc{u}) + \Div (\varrho \vc{u} \otimes \vc{u}) -\Grad p(\varrho)
= \mu\Delta\vc{u}+(\mu+\lambda)\Grad \Div(\vc{u}) +  \vc{H}\cdot\Grad \vc{H}-\frac{1}{2}\Grad(|\vc{H}|^2)+\varrho \vc{f},
\eF 
\bFormula{barop3}
\partial_t (\varrho \vartheta) + \Div (\varrho \vc{u} \vartheta) -\kappa\Delta\vartheta+p\Div(\vc{u})
= \frac{\mu}{2}|\Grad \vc{u}+\Grad^t\vc{u}|=\mu(\Div\vc{u})^2+\sigma(\Curl \vc{H})^2 +  \vc{H}\cdot\Grad \vc{H}-\frac{1}{2}\Grad(|\vc{H}|^2)+\varrho \vc{f},
\eF 
\bFormula{barop4}
\partial_t \vc{H} + \Div(\vc{u})\ \vc{H}+  \vc{u}\cdot \Grad \vc{H} ) =\nu \Delta \vc{H},\ \ \ \ \Div \vc{H} = 0,
\eF
\bFormula{barop5}
 (\varrho, \vc{u},\vc{H})_{t=0} = (\varrho_0, \vc{u}_0,\vc{H}_0),\ \ \ \ \Div \vc{H}_0 = 0.
\eF
Bian and Guo proved the following local existence result \cite{BG2} (see also \cite{BG} \cite{BY})
\bTheorem{BiGu}
Let $d\geq 2$, $\overline{\rho}>0$, $2\mu+d\lambda>0$, $\nu>0$.

Assume that the initial data satisfy $\varrho_0-\overline{\rho}\in  B^{d/p}_{p,1}$,
 $ \vc{u}_0\in (B^{d/p}_{p,1})^d$,
 $\vartheta_0\in B^{d/p}_{p,1}$ and
 $\vc{H}_0\in (B^{d/p}_{p,1})^d$,
with  $\varrho_0-\overline{\rho}$ bounded away from $0$.

Then there exists a positive number $T$ such that the system \eqref{barop1}-\eqref{barop4} has a unique solution $(\varrho,\vc{u},\vartheta,\vc{H}$ 
on the time interval $(0,T)$ such that
\[
\varrho-\overline{\rho}\in C([0,T);  B^{d/p}_{p,1}),
\]
\[
\vc{u}\in  \left(C([0,T); B^{d/p-1}_{p,1})\cap L^1([0,T); B^{d/p+1}_{p,1})\right)^d,
\]
\[
\vartheta\in  C([0,T); B^{d/p-2}_{p,1})\cap L^1([0,T); B^{d/p}_{p,1}),
\]
\[
\vc{H}\in  \left(C([0,T); B^{d/p-1}_{p,1})\cap L^1([0,T); B^{d/p+1}_{p,1})\right)^d.
\]
Moreover if $p\in[2,d]$ the solution is unique.
\eT



The first step of the proofs of both of these critical results is to use a smooth approximation of the data. For this regularized problem, previous
results of R. Danchin \cite{D2} allow to prove local existence using a linearization argument. Then, uniform bounds
are proved on the solution, allowing to pass to the limit as the regularization parameter vanishes. 

\section{One-dimensional models}
\label{1D}

In this section, we give an account of results in the
one-dimensional case, which allows for drastic simplification. In such a case, the use of Lagrangian variables
gives a much simpler system which is amenable by standard techniques. 

\bigskip

In dimension one, it is commonly assumed that the unknowns depend only on the first coordinate of
$x_1$ of $x$. The velocity of the fluid is still a three-dimensional vector,
but its first component plays a special role. We denote it by $u$, while the transverse velocity is $\vc w$:
\begin{displaymath}
  \vc u =
  \begin{pmatrix}
    u \\ w_1 \\ w_2 
  \end{pmatrix} =
  \begin{pmatrix}
    u \\ \vc w
  \end{pmatrix}.
\end{displaymath}

 Hence,
denoting $\vc B = (b_1,b_2,b_3)^T$, the constraint $\operatorname{div}(\vc B) = 0$ gives $\partial_1 b_1 = 1$, so
that $b_1$ is independent of $x_1$. Further, if $f$ is a vector field which is independent of $x_1$, then the first
component of $\Curl(f)$ vanishes. Hence, \eqref{b7} implies that $\partial_t b_1 = 1$. As a consequence, we may
assume without loss of generality that $b_1 = 1.$ Hence, denoting by $\vc b = (b_2,b_3)^T$ the transverse magnetic
field, that is,
\begin{displaymath}
  \vc B =
  \begin{pmatrix}
    b_1 \\ b_2 \\ b_3
  \end{pmatrix} =
  \begin{pmatrix}
    1 \\ \vc b
  \end{pmatrix},
\end{displaymath}
we find the following system
\begin{equation}\label{eq:1D}
\left\{
  \begin{aligned}
   & \partial_t \varrho + \partial_{x_1} \left(\varrho u\right) = 0, \\
   & \partial_t \left(\varrho u \right) + \partial_{x_1}\left(\varrho u^2 + p \right) =
    \partial_{x_1}\left[\left(\frac43\nu + \eta\right) \partial_{x_1}u\right] - \vc b \cdot \partial_{x_1}\left(\frac{\vc b}\mu\right), \\
    & \partial_t\left(\varrho\vc w \right) + \partial_{x_1} \left(\varrho u \vc w - \frac{\vc b}\mu\right) = \partial_{x_1}
    \left(\nu\partial_{x_1} \vc w \right), \\
    & \partial_t {\cal E} + \partial_{x_1} \left[\left({\cal E} + \frac {|\vc b|^2}{2\mu} + p \right)u - {\vc
        w}\cdot {\vc b}\right] = \partial_{x_1}\left[\left(\frac43\nu +
        \eta\right) u \partial_{x_1}u\right] + \partial_{x_1}\left( \kappa_F\partial_{x_1}\vartheta\right)
    \\ & \hspace{1cm}+ \partial_{x_1} \left(\frac{\vc b}\sigma\cdot \partial_{x_1}\left(\frac{\vc b}{\mu} \right) \right)
    + \partial_{x_1}\left(\nu {\vc w}\cdot\partial_{x_1} {\vc w}\right),\\
    & \partial_t\vc b + \partial_{x_1}\left(u\vc b - \vc w \right) = \partial_{x_1}\left(\frac 1
      \sigma \partial_{x_1}\frac{\vc b}\mu\right), 
  \end{aligned}\right.
\end{equation}
where the total energy $\cal E$ is given by 
\begin{displaymath}
  {\cal E} = \frac12 \varrho \left(u^2 + |\vc w|^2\right) + \varrho e + {\cal M}\left(\frac{|\vc b|}{\mu} \right).
\end{displaymath}
Recall that ${\cal M}$ is defined by \eqref{EM}, and that, if $\mu$ is a constant, we have 
\begin{displaymath}
 \displaystyle{\cal
  M}\left(\frac{|\vc b|}{\mu}\right) =\frac{|\vc b|^2}{2\mu}. 
\end{displaymath}
In what follows, we will always assume that $\mu$ is a constant. 

In the special case of a barotropic fluid, that is, without the energy equation (fourth equation of \eqref{eq:1D}
and with a pressure law given by $p(\varrho) = A \varrho^\gamma,$ $\gamma \geq 1$, it is proved in \cite{Z} that
the above system has a unique global smooth solution in the spatial domain $\Omega =(0,1)$, for initial data
$(\varrho_0, u_0,\vc w_0, \vc b_0)$ such that $\varrho_0\geq 0,$ $\varrho_0\in H^1$, $(u_0,\vc w_0)\in H^1_0\cap
H^2$, $\vc b_0 \in H^1_0$. Further, a compatibility condition is needed, in the spirit of those used in the case of
Navier-Stokes equations \cite{CK1,CK2,LXY}. All the coefficients $\nu,\eta,\sigma,\mu$ are assumed to be
constant. The method of proof is here again to prove local existence using Banach fixed-point theorem. Then a
priori estimates are proved, which allow to apply local existence again and extend the time interval to
infinity. These a priori estimates give, in a addition, a stability of the solution in Sobolev spaces. 

Next, in the paper \cite{ZX}, Zhang and Xie study the full system \eqref{eq:1D}, with additional gravitation terms:
\begin{equation}\label{eq:1D_gravite}
\left\{
  \begin{aligned}
   & \partial_t \varrho + \partial_{x_1} \left(\varrho u\right) = 0, \\
   & \partial_t \left(\varrho u \right) + \partial_{x_1}\left(\varrho u^2 + p \right) =
    \partial_{x_1}\left[\left(\frac43\nu + \eta\right) \partial_{x_1}u\right] - \vc b
    \cdot \partial_{x_1}\left(\frac{\vc b}\mu\right) + \varrho\partial_{x_1}\psi, \\
    & \partial_t\left(\varrho\vc w \right) + \partial_{x_1} \left(\varrho u \vc w - \frac{\vc b}{\mu}\right) = \partial_{x_1}
    \left(\nu\partial_{x_1} \vc w \right), \\
    & \partial_t {\cal E} + \partial_{x_1} \left[\left({\cal E} + \frac {|\vc b|^2}{2\mu} + p \right)u - {\vc
        w}\cdot {\vc b}\right] = \partial_{x_1}\left[\left(\frac43\nu +
        \eta\right) u \partial_{x_1}u\right] + \partial_{x_1}\left( \kappa_F\partial_{x_1}\vartheta\right)
    \\ & \hspace{1cm}+ \partial_{x_1} \left(\frac{\vc b}\sigma\cdot \partial_{x_1}\left(\frac{\vc b}{\mu} \right) \right)
    + \partial_{x_1}\left(\nu {\vc w}\cdot\partial_{x_1} {\vc w}\right)+
    \varrho u \partial_{x_1}\psi ,\\
    & \partial_t\vc b + \partial_{x_1}\left(u\vc b - \vc w \right) = \partial_{x_1}\left(\frac 1
      \sigma \partial_{x_1}\frac{\vc b}\mu\right), 
  \end{aligned}\right.
\end{equation}
In this system, $\psi$ is the gravitation field, defined by Poisson's equation:
\begin{equation}
  \label{eq:poisson}
\left\{
  \begin{aligned}
    &-\partial_{x_1}^2 \psi = G \rho \text{ in } \Omega, \\
    & \psi = 0 \text{ on } \partial\Omega.
  \end{aligned}
\right.
\end{equation}
where $G>0$ is the gravitational constant. The conductivity $\kappa_F = \kappa_F(\varrho,\vartheta)$ is not assumed to
be a constant, but satisfies the growth assumptions
\begin{equation}
  \label{eq:borne_kappa}
  \frac 1 C\left(1+\vartheta^q\right)\leq \kappa_F, \quad \left|\frac{\partial \kappa_F}{\partial \varrho}\right| \leq 1+
  \vartheta^q, \quad \forall q>\frac52,
\end{equation}
for some constant $C>0$. Further, the equation of state includes the radiative term $\vartheta^4$, that is
\begin{equation}\label{eq:eos}
  p = R\rho\vartheta + \frac a 3 \vartheta^4, \quad e = C_V \vartheta + \frac a\varrho \vartheta^4.
\end{equation}
The method is similar to \cite{Z}. This result has been improved by Qin and Hu \cite{QH},
in which $q> \frac{2+\sqrt{211}}9 \approx 0.3$ is allowed.

\medskip

Another approach in this one-dimensional setting is to use Lagrangian coordinates. Here, we assume that we use a
bounded domain, which, without loss of generality, we may assume to be $(0,1)$. We also assume that $\int_0^1
\rho_0(x_1) dx_1 = 1.$ Setting
\begin{equation}
  \label{eq:y_lagrangien}
  y = y(t,x_1) = \int_0^{x_1} \varrho(t,\xi)d\xi,
\end{equation}
we then have $0\leq y \leq 1$, and, due to mass conservation
\begin{displaymath}
  \int_0^1 \varrho(t,x_1) dx_1 = \int_0^1 \varrho_0(t,x_1)dx_1=1.
\end{displaymath}
Thus, defining $v=1/\rho$ the specific volume of the fluid, we immediately have, using Dirichlet Boundary
conditions
\begin{displaymath}
  \partial_t y + \varrho u = 0.
\end{displaymath}
Hence, system \eqref{eq:1D} becomes
\begin{equation}
  \label{eq:1D_lagrangien}
  \left\{
    \begin{aligned}
      & \partial_t v - \partial_y u = 0, \\
      & \partial_t u + \partial_y p = \partial_y \left(\frac{\frac43\nu+\eta}v \partial_y u\right) - {\vc
        b}\cdot \partial_y \left(\frac{\vc b}{\mu}\right), \\
      & \partial_t {\vc w} - \partial_y \left(\frac{\vc b}\mu\right) = \partial_y \left(\frac \nu v \partial_y {\vc w}\right), \\
      & \partial_t \left( v{\cal E}\right) + \partial_y \left(pu + \frac{|{\vc b}|^2}{2\mu}u - {\vc w}\cdot {\vc b}
      \right) =  \partial_y\left[\left(\frac43 \nu + \eta\right) u \frac{\partial_y u}v + \kappa_F \frac{\partial_y
          \vartheta}v  + \frac{{\vc b}}{\sigma v} \partial_y \left(\frac{\vc b}{\mu}\right) + \frac{\nu}v \vc w
        \cdot \partial_y \vc w\right], \\
      & \partial_t \left(v\vc b \right) - \partial_y {\vc w} = \partial_y \left( \frac 1 {\sigma v} \partial_y \left(
          \frac{\vc b}{\mu} \right)\right).
    \end{aligned}
\right.
\end{equation}
In principle, this system is set on the interval $\Omega = (0,1)$ with homogeneous Dirichlet boundary conditions:
\begin{equation}
  \label{eq:bc_lagrange}
  \left(u,\vc w, \vc b, \partial_y \vartheta\right) = 0\text{ on }\partial\Omega.
\end{equation}
However, the corresponding Cauchy problem ($\Omega = \R$) is also studied by some authors, with the same kind of
boundary conditions. System \eqref{eq:1D_lagrangien}-\eqref{eq:bc_lagrange} is completed by initial conditions
$(v_0, u_0, \vc w_0, \vartheta_0, \vc b_0)$.

\medskip

The first existence result concerning the Cauchy problem \eqref{eq:1D_lagrangien} is the global existence of
Kawashima and Okada \cite{KO}. The equation of state is assumed to satisfy natural conditions similar to \eqref{b23}. It is proved
in \cite{KO} that, given a constant state $\left(\overline v, \overline u=0 , \overline{\vc w}=0, \overline\vartheta,
  \overline{\vc b}\right),$ (with $\overline v>0$ and $\overline \vartheta >0$) if the initial condition is close
to this state in $H^2(\R)$, then system \eqref{eq:1D_lagrangien}-\eqref{eq:bc_lagrange} has a unique global smooth
solution $\left(v,u,{\vc w},\vartheta, \vc b\right)$ such that, for any $t\geq 0,$ $\left(v(t) -\overline v, u(t)
  ,\vc w(t) - \overline{\vc w}, \vartheta(t) - \overline\vartheta, \vc b(t) -
  \overline{\vc b} \right) \in H^2(\R)$. Here again, local existence is proved using classical results, and global
energy estimates are proved for these smooth solutions. This allows to extend the existence result to any positive
time. The result may be summarized as follows:

\begin{Theorem}[Kawashima-Okada, 1982 \cite{KO}]
Consider system \eqref{eq:1D_lagrangien} set in $\Omega = \R$. Assume that the equation of state satisfies
\begin{displaymath}
  \frac{\partial p}{\partial \varrho} >0, \quad \frac{\partial e}{\partial \vartheta} >0,
\end{displaymath}
and that $de = \vartheta d S - p d\left(\frac 1 \varrho\right),$ where $S$ is the entropy of the fluid. Assume that
$\mu$ is a positive constant, and that $\eta,\nu,\kappa_F$ are positive smooth functions of $\varrho$ and
$\vartheta$. Let 
\begin{displaymath}
  \left(\overline v, \overline u, \overline{\vc w}, \overline\vartheta,\overline{\vc b}\right) = \left(\overline v, 0, 0, \overline\vartheta,0\right) 
\end{displaymath}
be a constant equilibrium state of \eqref{eq:1D_lagrangien}. Let $\left(v_0, u_0, {\vc w}_0, \vartheta_0, {\vc
    b}_0\right)$ be an initial data such that
  \begin{displaymath}
    \left(v_0 - \overline v, u_0, {\vc w}_0,  \vartheta_0 -\overline\vartheta,{\vc b}_0\right)\in H^2(\R).
  \end{displaymath}
If $\left\|v_0 - \overline v, u_0, {\vc w}_0,  \vartheta_0 -\overline\vartheta,{\vc b}_0 \right\|_{H^2(\R)}$ is sufficiently small, then \eqref{eq:1D_lagrangien} has a unique global smooth
solution. 
\end{Theorem}

Actually, coefficients are allowed to vanish, provided they are constant: the situation allowed are (a)
$\nu=\eta=0$, $\mu,\kappa_F>0$, (b) $\nu,\eta,\mu>0$, $\kappa_F = 0$, (c) $\eta=\nu=\kappa_F = 0$,
$\mu>0$.

Another result can be found in \cite{CW}. In this contribution, system \eqref{eq:1D} is considered on the domain
$\Omega= \Omega(t)= (0,x(t))$, where the free boundary $x(t)$ is defined by $x'(t) = u(t,x(t))$. The boundary
conditions are given by
\begin{equation}
  \label{eq:bc_Chen_Wang}
  \left(\vc w, \vc b, \partial_{x_1}\vartheta\right)_{\partial\Omega(t)} = 0, \quad u(t,0) = 0, \quad \left[p -
      \left(\frac43\nu+\eta \right) \partial_{x_1} u \right]_{x_1 = x(t)} = 1.
\end{equation}
Then, the Lagrangian coordinates are used to transform the system into \eqref{eq:1D_lagrangien}, in the domain
$\Omega = (0,1)$, with the boundary conditions 
\begin{equation}
  \label{eq:bc_Chen_Wang_lagrange}
  \left(\vc w, \vc b, \partial_{y}\vartheta\right)_{\partial\Omega} = 0, \quad u(t,0) = 0, \quad \left[p -
      \frac 1 v \left(\frac43\nu+\eta \right) \partial_{y} u \right]_{y=1} = 1.
\end{equation}
We are back to a system which is similar to \eqref{eq:1D_lagrangien}-\eqref{eq:bc_lagrange}, to which the strategy
of Kawashima-Okada \cite{KO} is applied, thereby proving global existence if the initial data is close to a given
constant state, in Sobolev spaces. Moreover, in the special case of a perfect gas equation of state, a global
existence of weak solutions is proved. For this purpose, the initial data, which may be discontinuous, is
approximated by smooth functions. This allows to apply the global smooth existence result. The estimates
established in the smooth case are then used again to prove that one can pass to the limit as the regularization
parameter tends to $0$. This is based on Aubin-Lions Lemma (see \cite{AU,LI}). All these results are proved using
the assumption that $\sigma=1$, $\eta,\nu$ are bounded smooth functions of $v$, that the energy and pressure of
the fluid satisfy some polynomial growth property as functions of $\vartheta$, and the fact that 
\begin{displaymath}
  \kappa_F(v,\vartheta) + \left| \partial_v \kappa_F(v,\vartheta)\right| + \left| \partial_v^2
    \kappa_F(v,\vartheta)\right| \leq C \left(1+\vartheta^{q}\right),
\end{displaymath}
for some constant $C>0$ and some power $q\geq 2.$ This is in sharp contrast with \cite{KO}, in which $\kappa_F$ is
assumed to be constant. 

Using a similar method, Wang \cite{W} proved that system~\eqref{eq:1D_lagrangien}-\eqref{eq:bc_lagrange} has a
global weak solution in $H^1(\Omega)$. The same kind of assumptions are used on the equation of state and the
diffusion coefficient $\kappa_F$. 

In a similar setting, it is proved in \cite{CW2} that the solutions depend continuously on the initial data. 

\medskip

In the case of constant coefficients, with an equation of state satisfying 
\begin{displaymath}
  \left|\frac{\partial
    p}{\partial e}(v,e)\right| + e(\vartheta) + \left|\frac{\partial e}{\partial \vartheta}(v,\vartheta)\right|\leq C,
\end{displaymath}
the Cauchy problem \eqref{eq:1D_lagrangien} set on $\Omega = \R$ is studied by Hoff and Tsyganov in \cite{HT,HT2}. It is proved that weak solutions exist and are unique, provided the total energy is small. Under
an additional integrability assumption on $\partial_t \vc b$, continuity with respect to the initial data is proved
(in $L^2([0,T]\times \R)$, for any $T>0$).

\medskip

In \cite{FJN}, for system \eqref{eq:1D_lagrangien}-\eqref{eq:bc_lagrange}, the limit of vanishing shear viscosity
$\nu\to 0$ is studied, under the same assumptions as in \cite{CW}. It is proved that the solution converges to that
of the limit problem strongly in $L^2([0,T]\times \Omega)$ for $\rho,\vc w,\vartheta$, and strongly in $L^2([0,T],
H^1(\Omega))$ for $u,\vc b$.

\medskip

Finally, still in Lagrangian coordinates, but in the special case $\vc w = 0$,
system~\eqref{eq:1D_lagrangien}-\eqref{eq:bc_lagrange} has been studied by Iskenderova \cite{I}. Assuming constant
coefficients and a perfect gas equation of state, Iskenderova proves that the system has a unique weak
solution. Note that the initial density $\rho_0$ is allowed to vanish, provided that $\partial_y
\left(\rho_0^{1/4}\right) \in L^2$, and that
\begin{displaymath}
  \rho_0(y) \leq C\inf_{0\leq z \leq y} \rho_0(z),
\end{displaymath}
for some constant $C>0$. Here again, the method of proof is a regularization of the system: existence and
uniqueness is easily proved for a system in which $\rho_0\geq \varepsilon>0$. Then, the limit $\varepsilon \to 0$
is studied using a priori estimates.

\bibliographystyle{spbasic}

\end{document}